\documentclass[11pt]{article}
\usepackage{lmodern}
\usepackage[T1]{fontenc}
\usepackage{CJKutf8}

\usepackage{graphicx}
\usepackage{float}
\usepackage{amsmath,amsthm,amsfonts,amssymb}
\usepackage{ parskip }
\usepackage{ authblk }
\usepackage{multicol}
\usepackage{caption}
\usepackage{subcaption}
\usepackage[colorlinks=true,linkcolor=blue,citecolor=blue,urlcolor=blue,bookmarks=true]{hyperref}
\usepackage[margin=1.1cm, top=1.5cm, bottom=1.5cm]{geometry}
\usepackage[utf8 ]{inputenc}
\numberwithin{equation}{section}

\begin{document}
	\lineskip 1.5em
	\title{ Love Dynamical Model with persepectives of Piecewise Differential Operators}
	\author[]{Atul Kumar } 
	\affil[]{Dayalbagh Educational Institute, Department of Mathematics, India}
	\date{~}
	\maketitle
	\begin{abstract}
		For love dynamical models, a new idea combining piecewise concept for intger-order, stochastic, and fractional derivatives is presented in order to capture the chaos and several crossover emotional scenerios. Under the assumptions of linear growth and Lipschitz condition, the fixed-point theorem explain the uniqueness and existence to the models under the investigation. The piecewise derivatives were approximated utilising the Lagrange interpolation method, and the computer results were demonstrated numerically for several values of order $\alpha$. It was observed that the recently presented new idea in love dynamical models can represent disordered emotional patterns in passionate loving partnerships. 
	\end{abstract}

	\section{Introduction}
	Mathematical, psychological, and biological scientists have worked together to develop mathematical models that illustrate the crossover emotional patterns in romantic and interpersonal relationships. Using integer-order derivatives, the majority of dynamic model formulations have been presented \cite{ Barkley, Bielczyk, Rinaldi2, dengLiao}. As an illustration of the crossover emotional behaviors between two partners under the influence of family on each individual, Liu et al. \cite{Liu} suggested a novel love dynamical model. Recently, Yifan Liu and colleagues \cite{yifanliu} created a new nonlinear love model with two delays by injecting nonlinearity originating from a logistic component. The linearization theory was employed by the authors to investigate its stability and Hopf bifurcation. For more noteworthy results, the reader can consult the references mentioned above. 
	
	Fractional calculus has become a more effective tool than classical calculus in the last few decades for the purpose of explaining and illuminating nonlinear phenomena that arise in a variety of engineering and applied science fields \cite{Din, Babaei, Zhang, Danane}. Among them, it is important to note the relevance of the fractional approach to the analysis of love-dynamical models \cite{song1}, and the recently established formulation may represent peculiar emotional patterns. However, the task of characterizing crossover behaviors remains difficult. The kernels of several fractional derivatives and integrals were examined \cite{Khan, Baleanu, Peter, Atangana3}.
	For illustrating crossover behaviors over the power-law kernel, they discovered that the generalized Mittag-Leffler function and exponential decay are more appropriate. Nonetheless, considerable constraints are necessary to illustrate the transition from a deterministic to a stochastic context \cite{Danane2}. 
	
	The intricacy of nature often presents us with multifaceted real-world situations that seem beyond the scope of current theories and tools. One such example is a phenomenon that shows behavior that shifts from fuzzy to stochastic or from stochastic to power law. For that reason, Atangana et al. \cite{Atangana4} presented a novel idea of piecewise derivatives and integrals and shown how effective it is in representing crossover in models that are chaotic and epidemiological. The Covid-19 model in the third wave \cite{Atangana3} and rhythmic heartbeats \cite{Atangana6} are some examples of the piecewise pattern that the piecewise modeling idea is explored to characterize the crossover behavior. Author applied this idea to capture chaos and crossover emotional patterns in love dynamical systems due to its effectiveness and adaptability. Piecewise concept-based crossover behaviors for a model of tumor growth and radiation response \cite{Arık}. COVID-19 infection dynamics with application of piecewise fractional differential equation \cite{Li XP}. Piecewise fractional examination of the effect of migration in interactions between plants, pathogens, and herbivores \cite{ur Rahman}. Analysis of Ulam-Hyers stability and solvability for nonlinear piecewise fractional cancer dynamic systems \cite{KhanS} (). A case study in Yemen: examination and modeling of the piecewise dynamics model of malaria transmission ~\cite{Aldwoah}.

	Following is how author constituted the study. First, author reviewed the background information and the meanings of integral and piecewise differential operators in Section 2. Section 3 presents the modified love dynamical models in piecewise sense. In Section 4, author prove results on the uniqueness and existence for the system. In Section 5, numerical solutions and computer results are shown for the considered models at different values of arbitrary fractional order $\alpha$. The conclusion follows Section 6.

	\section{Basic Definitions}
	\textbf{Definition 1:} 
	Let $\alpha \in (0,1)$ be a real number and $U$ be a differentiable function, then the Caputo
	derivative\cite{Owolabi} with order $\alpha$ and lower-limit $a$ is defined as
	\begin{equation*}	
		_a^CD^{\alpha}_tU(t) = 
		\frac{1}{\Gamma{(1-\alpha)}}\int_a^t(t-s)^{-\alpha}  { U'(s)} d{s} 
	\end{equation*}
	
	\textbf{Definition 2:} Assume $U \in H^{1}(u_{1}, u_{2})$, for $u_{2}>u_{1}$ and  $\alpha \in (0, 1]$ with $U$ being a function then, the Atangana–Baleanu fractional derivative\cite{Owolabi} with arbitrary order $\alpha$ and lower-limit $a$ in the sense of Caputo is defined as 
	\begin{equation*}	
		_a^{ABC}D^{\alpha}_tU(t) = 
		\frac{AB(\alpha)}{1-\alpha}\int_a^t U'(s) E_{\alpha}\left[-\frac{\alpha}{1-\alpha}(t-s)^{\alpha}\right] ds
	\end{equation*}
	where $AB(\alpha)=1-\alpha+\frac{\alpha}{\Gamma(\alpha)}$ and $E_{\alpha}(z)=\displaystyle\sum_{n=0}^{\infty}\frac{z^n}{\Gamma(\alpha n + 1)}$ are known as the normalization function and the Mittag-Leffler function, respectively. 
	
	\textbf{Definition 3:} Let $U \in H^{1}(u_{1}, u_{2})$, for $u_{2}>u_{1}$   and $\alpha \in (0, 1]$ then the Caputo–Fabrizio fractional derivative\cite{Owolabi} of arbitrary order $\alpha$ with exponential law and lower-limit $a$ is defined by
	\begin{equation*}	
		_a^{CF}D^{\alpha}_tU(t) = 
		\frac{M(\alpha)}{1-\alpha}\int_a^t U'(s) \exp\left[-\frac{\alpha}{1-\alpha}(t-s)\right] ds
	\end{equation*}
	where $M(\alpha)$ is a normalization function such that M(0)=M(1)=1
	
	\textbf{Definition 4:} For a function $U$, a piecewise integral is given as \cite{Atangana4} 
	\begin{equation}
		_0^{PPL}J^{\alpha}_tU(t) = 
		\begin{cases}
			\displaystyle\int_0^{t} U(s) d{s} & 0\leq t\leq t_{1} \\ 
			\frac{1}{\Gamma{(\alpha)}}\displaystyle\int_{t_{1}}^t U(s)(t-s)^{\alpha-1}    d{\tau}	& t_{1}\leq t\leq a
		\end{cases}
	\end{equation}
	
	\textbf{Definition 5: } 
	For a function $U$, another definition of a piecewise integral is given as \cite{Atangana4}  
	\begin{equation}
		_0^{PCF}J^{\alpha}_tU(t) = 
		\begin{cases}
			\displaystyle\int_0^{t} U(s) d{s} & 0\leq t\leq t_{1} \\ 
			\frac{1-\alpha}{M(\alpha)}U(t)+\frac{\alpha}{M(\alpha)}\displaystyle\int_{t_{1}}^t U(s)    d{s}	& t_{1}\leq t\leq a
		\end{cases}
	\end{equation}
	
	\textbf{Definition 6: } 
	For a function $U$, another definition of a piecewise integral is given as \cite{Atangana4}  
	\begin{equation}
		_0^{PAB}J^{\alpha}_tU(t) = 
		\begin{cases}
			\displaystyle\int_0^{t} U(s) d{s} & 0\leq t\leq t_{1} \\ 
			\frac{1-\alpha}{AB(\alpha)}U(t)+\frac{\alpha}{AB(\alpha)}\frac{1}{\Gamma{(\alpha)}}\displaystyle\int_{t_{1}}^t U(s)(t-s)^{\alpha-1}    d{s}	& t_{1}\leq t\leq a
		\end{cases}
	\end{equation}
	
	\textbf{Definition 7: }For a function $U$, a piecewise derivative with power-law kernel will be defined as \cite{Atangana4}  
	\begin{equation}
		_0^{PC}D^{\alpha}_tU(t) = 
		\begin{cases}
			U'(t) & 0\leq t\leq t_{1} \\
			_{t_{1}}^{C}D^{\alpha}_tU(t)
			& t_{1}\leq t\leq a
		\end{cases}
	\end{equation}

	\textbf{Definition 8: }For a function $U$, a piecewise derivative with exponential-decay kernel is given as \cite{Atangana4} 
	\begin{equation}
		_0^{PCF}D^{\alpha}_tU(t) = 
		\begin{cases}
			U'(t) & 0\leq t\leq t_{1} \\
			_{t_{1}}^{CF}D^{\alpha}_tU(t)
			& t_{1}\leq t\leq a
		\end{cases}
	\end{equation}
	
	\textbf{Definition 9:} 
	For a function $U$, a piecewise derivative with Mittag-Leffler kernel is given as \cite{Atangana4} 
	\begin{equation}
		_0^{PAB}D^{\alpha}_tU(t) = 
		\begin{cases}
			U'(t) & 0\leq t\leq t_{1} \\
			_{t_{1}}^{ABC}D^{\alpha}_tU(t)
			& t_{1}\leq t\leq a
		\end{cases}
	\end{equation}

	\section{Modeling of the love dynamical system}
	This section introduces an alternative approach to modeling a love dynamical system that displays feelings of romantic love relationships with diverse patterns. Assume, for the sake of argument, that a dynamical model with $k$ classes explains many patterns of love relationships without losing generality. Assume, without sacrificing generality, that distinct patterns will be taken into consideration for randomness, nonlocal processes, classical mechanical processes, and their permutations. The selection of piecewise differential operators in a real-world scenario where various patterns are seen will be guided by the collected data.
	
	\subsection*{Case 1: Classical-power-law-randomness\cite{Atangana3}}
	Here, we take it that romantic partnerships with piecewise love dynamics follow a classical-power-law-randomness process. Assume that the piecewise differential operator is injected into the first process, which goes from 0 to $a_{1}$. The second is the Caputo derivative, which has an origin at $a_{1}$ and a destination at $a_{2}$. Lastly, the stochastic differential equation that starts at $a_{2}$ and finishes at $a$ is a part of the third process. The complete process' mathematical expression for this instance will be provided as
	\begin{align}\label{linear_caputo}
		\begin{cases}
			&U'(t) = e(t,U_{i}), 0 \le t\le a_{1} ,  \\
			& U_{i}(0)=U_{i,0}, i=1,2,...,k  \\
			& _{a_{1}}^CD^{\alpha}_tU_{i}(t) = e(t,U_{i}), a_{1} \le t\le a_{2}  \\
			& U_{i}(a_{1})=U_{i,1}, i=1,2...k  \\
			&dU_{i}(t)=e(t,U_{i})dt+\sigma_{i}U_{i}d \mathbb{B}_{i}(t), a_{2} \le t\le a  \\
			& U_{i}(a_{2})=U_{i,2}, i=1,2...k  \\
		\end{cases}
	\end{align}
	where $\mathbb{B}_{i}(t)$, $i=1,2$, denote the standard Brownian Motions defined on a complete probability space $\left(\Omega, \mathcal{F}, \left\{\mathcal{F}_t\right\}_{t\geq 0}, P\right)$ with the filtration $\left\{\mathcal{F}_t\right\}_{t\geq 0}$ satisfying the usual conditions. Also, $\sigma_{i}$ denotes the intensity of $\mathbb{B}_{i}(t)$.
	\subsection*{Case 2:Classical-Mittag-Leffler-law-randomness\cite{Atangana3}}
	In this instance, we considered that a classical-Mittag-Leffler-law-randomness process underlies the romantic interactions in piecewise love dynamics. The second step differs from Case 1 in that the power-law kernel is not injected; instead, the Mittag-Leffler kernel is. The Atangana-Balaneu derivative takes the role of the Caputo derivative in the Caputo meaning. The following differential equations govern the entire process in the current case:
	
	\begin{align}\label{linear_caputo_mittag}
		\begin{cases}
	&U'(t) = e(t,U_{i}), 0 \le t\le a_{1} ,  \\
	& U_{i}(0)=U_{i,0}, i=1,2,...,k  \\
	& _{a_{1}}^{ABC}D^{\alpha}_tU_{i}(t) = e(t,U_{i}), a_{1} \le t\le a_{2}  \\
	& U_{i}(a_{1})=U_{i,1}, i=1,2...k  \\
	&dU_{i}(t)=e(t,U_{i})dt+\sigma_{i}U_{i}d \mathbb{B}_{i}(t), a_{2} \le t\le a  \\
	& U_{i}(a_{2})=U_{i,2}, i=1,2...k  \\
\end{cases}
	\end{align}
	
	\subsection*{Case 3:Classical-fading memory-randomness\cite{Atangana3}}
	In this scenario, we suppose that fading memory randomization processes are present in romantic relationships involving piecewise love dynamics. The second process, which has an exponential-decay kernel, differs from the situations mentioned previously in that we took into account the Caputo Fabrizio derivative. The following differential equations control the entire procedure in the current case:
	
	\begin{align}\label{linear_caputo_fading}
		\begin{cases}
	&U'(t) = e(t,U_{i}), 0 \le t\le a_{1} ,  \\
	& U_{i}(0)=U_{i,0}, i=1,2,...,k  \\
	& _{a_{1}}^{CF}D^{\alpha}_tU_{i}(t) = e(t,U_{i}), a_{1} \le t\le a_{2}  \\
	& U_{i}(a_{1})=U_{i,1}, i=1,2...k  \\
	&dU_{i}(t)=e(t,U_{i})dt+\sigma_{i}U_{i}d \mathbb{B}_{i}(t), a_{2} \le t\le a  \\
	& U_{i}(a_{2})=U_{i,2}, i=1,2...k  \\
\end{cases}
	\end{align}

	\section*{Mathematical Models}
	We have worked with two love dynamical models in this section: one for a linear situation and one for a nonlinear one. For Cases 1, 2, and 3, these mathematical models have been examined.
	\subsection{First Mathematical Model\cite{Owolabi}}
	The two species in the love dynamical model that follows are $r$ and $s$. When $t$ is reached, $r(t)$ and $s(t)$ denote the love or dislike for specie-$s$ and specie-$r$, respectively; $\rho_i$ for $i = 1, 2$ represents the forgetting constant; $\psi_i$ and $\gamma_i$ for $i = 1, 2$ represent the pleas for the love towards species $r$ and $s$, respectively; $\omega_1,$ $\omega_2$ represent the reactance for the love towards $r$ and $s$, respectively.
	\begin{align}\label{caputolinear}
		\begin{cases}
			& \frac{dr(t)}{dt} = -\rho_{1}r+\omega_{1}s+\gamma_{1} \psi_{2} \\
			& \frac{ds(t)}{dt} = -\rho_{2}s+\omega_{2}r+\gamma_{2} \psi_{1} \\
			& _{a_{1}}^CD^{\alpha}_tr(t) = -\rho_{1}r(t)+\omega_{1}s(t)+\gamma_{1} \psi_{2} \\
			&  _{a_{1}}^CD^{\alpha}_ts(t)= -\rho_{2}s+\omega_{2}r+\gamma_{2} \psi_{1} \\
			& dr(t)=\bigg(-\rho_{1}r+\omega_{1}s+\gamma_{1} \psi_{2}\bigg)dt+\sigma_{1}r(t)d \mathbb{B}_{1}(t) \\
			& ds(t)=\bigg(-\rho_{2}s+\omega_{2}r+\gamma_{2} \psi_{1}\bigg)dt+\sigma_{2}s(t)d \mathbb{B}_{2}(t) 
		\end{cases}
	\end{align}
	\begin{align}\label{ABlinear}
		\begin{cases}
	& \frac{dr(t)}{dt} = -\rho_{1}r+\omega_{1}s+\gamma_{1} \psi_{2} \\
	& \frac{ds(t)}{dt} = -\rho_{2}s+\omega_{2}r+\gamma_{2} \psi_{1} \\
	& _{a_{1}}^{ABC}D^{\alpha}_tr(t) = -\rho_{1}r(t)+\omega_{1}s(t)+\gamma_{1} \psi_{2} \\
	&  _{a_{1}}^{ABC}D^{\alpha}_ts(t)= -\rho_{2}s+\omega_{2}r+\gamma_{2} \psi_{1} \\
	& dr(t)=\bigg(-\rho_{1}r+\omega_{1}s+\gamma_{1} \psi_{2}\bigg)dt+\sigma_{1}r(t)d \mathbb{B}_{1}(t) \\
	& ds(t)=\bigg(-\rho_{2}s+\omega_{2}r+\gamma_{2} \psi_{1}\bigg)dt+\sigma_{2}s(t)d \mathbb{B}_{2}(t) 
\end{cases}
	\end{align}
	\begin{align}\label{CFlinear}
		\begin{cases}
	& \frac{dr(t)}{dt} = -\rho_{1}r+\omega_{1}s+\gamma_{1} \psi_{2} \\
	& \frac{ds(t)}{dt} = -\rho_{2}s+\omega_{2}r+\gamma_{2} \psi_{1} \\
	& _{a_{1}}^{CF}D^{\alpha}_tr(t) = -\rho_{1}r(t)+\omega_{1}s(t)+\gamma_{1} \psi_{2} \\
	&  _{a_{1}}^{CF}D^{\alpha}_ts(t)= -\rho_{2}s+\omega_{2}r+\gamma_{2} \psi_{1} \\
	& dr(t)=\bigg(-\rho_{1}r+\omega_{1}s+\gamma_{1} \psi_{2}\bigg)dt+\sigma_{1}r(t)d \mathbb{B}_{1}(t) \\
	& ds(t)=\bigg(-\rho_{2}s+\omega_{2}r+\gamma_{2} \psi_{1}\bigg)dt+\sigma_{2}s(t)d \mathbb{B}_{2}(t) 
\end{cases}
	\end{align}
	
	\subsection{Second Mathematical Model \cite{Owolabi}}
	The term involving $\epsilon$ captures the compensations for conflict with pity or flattery in this mathematical model, which is the nonlinear variant of the linear model. Similar emotional patterns to the preceding situation are described by additional notations.
	\begin{align}\label{caputononlinear}
		\begin{cases}
			& \frac{dr(t)}{dt} = -\rho_{1}r+\omega_{1}s(1-\epsilon s^2)+\psi_{1} \\
			&	 \frac{ds(t)}{dt} = -\rho_{2}s+\omega_{2}r(1-\epsilon r^2)+\psi_{2} \\
			& _{a_{1}}^CD^{\alpha}_tr(t) = -\rho_{1}r+\omega_{1}s(1-\epsilon s^2)+\psi_{1} \\
			& _{a_{1}}^CD^{\alpha}_ts(t) = -\rho_{2}s+\omega_{2}r(1-\epsilon r^2)+\psi_{2} \\
			& 		dr(t)=\bigg(-\rho_{1}r+\omega_{1}s(1-\epsilon s^2)+\psi_{1} \bigg)dt+\sigma_{1}r(t)d \mathbb{B}_{1}(t) \\
			&		ds(t)=\bigg( -\rho_{2}s+\omega_{2}r(1-\epsilon r^2)+\psi_{2}\bigg)dt+\sigma_{2}s(t)d \mathbb{B}_{2}(t) 
		\end{cases}
	\end{align}
	\begin{align}\label{ABnonlinear}
		\begin{cases}
	& \frac{dr(t)}{dt} = -\rho_{1}r+\omega_{1}s(1-\epsilon s^2)+\psi_{1} \\
	&	 \frac{ds(t)}{dt} = -\rho_{2}s+\omega_{2}r(1-\epsilon r^2)+\psi_{2} \\
	& _{a_{1}}^{ABC}D^{\alpha}_tr(t) = -\rho_{1}r+\omega_{1}s(1-\epsilon s^2)+\psi_{1} \\
	& _{a_{1}}^{ABC}D^{\alpha}_ts(t) = -\rho_{2}s+\omega_{2}r(1-\epsilon r^2)+\psi_{2} \\
	& 		dr(t)=\bigg(-\rho_{1}r+\omega_{1}s(1-\epsilon s^2)+\psi_{1} \bigg)dt+\sigma_{1}r(t)d \mathbb{B}_{1}(t) \\
	&		ds(t)=\bigg( -\rho_{2}s+\omega_{2}r(1-\epsilon r^2)+\psi_{2}\bigg)dt+\sigma_{2}s(t)d \mathbb{B}_{2}(t) 
\end{cases}
	\end{align}
	\begin{align}\label{CFnonlinear}
		\begin{cases}
	& \frac{dr(t)}{dt} = -\rho_{1}r+\omega_{1}s(1-\epsilon s^2)+\psi_{1} \\
	&	 \frac{ds(t)}{dt} = -\rho_{2}s+\omega_{2}r(1-\epsilon r^2)+\psi_{2} \\
	& _{a_{1}}^{CF}D^{\alpha}_tr(t) = -\rho_{1}r+\omega_{1}s(1-\epsilon s^2)+\psi_{1} \\
	& _{a_{1}}^{CF}D^{\alpha}_ts(t) = -\rho_{2}s+\omega_{2}r(1-\epsilon r^2)+\psi_{2} \\
	& 		dr(t)=\bigg(-\rho_{1}r+\omega_{1}s(1-\epsilon s^2)+\psi_{1} \bigg)dt+\sigma_{1}r(t)d \mathbb{B}_{1}(t) \\
	&		ds(t)=\bigg( -\rho_{2}s+\omega_{2}r(1-\epsilon r^2)+\psi_{2}\bigg)dt+\sigma_{2}s(t)d \mathbb{B}_{2}(t) 
\end{cases}
	\end{align}

	\section{Existence and uniqueness of Love Dynamical model\cite{Atangana3}}
	Specifically, we prove the existence and uniqueness of the system solution of the linear case with $k=2$, which is a special instance of the first piecewise love dynamical model \eqref{linear_caputo}. We saw that different strategies are needed for each interval. Thus, for the first and second situations, we used the Banach fixed-point theorem to apply the well-known approach. On the other hand, a different method involving the stochastic differential equation is used in the last instance. To keep things simple, we demonstrate the existence and uniqueness of the subsequent piecewise mathematical model. 
	\begin{align*}
		\begin{cases}
			& U'(t) = e(t,U),   0 \le t\le a_{1}\\
			& _{a_{1}}^CD^{\alpha}_tU(t) =e(t,U),  a_{1} \le t\le a_{2}\\
			& dU(t)=e(t,U)dt+\sigma_{i}Ud \mathbb{B}(t), a_{2} \le t\le a
		\end{cases}
	\end{align*}
	where	$\sigma=(\sigma_{1}(t),\sigma_{2}(t))$, $\mathbb{B}=(\mathbb{B}_{1}(t),\mathbb{B}_{2}(t))$, $U=(r,s)$, $e=(e_{1},e_{2})$, $e_{1}=-\rho_{1}r+\omega_{1}s(1-\epsilon s^2)+\psi_{1}$, and 
	$ e_{2} = -\rho_{2}s+\omega_{2}r(1-\epsilon r^2)+\psi_{2} $. Our findings demonstrate that $e_{i}(t,r,s)$ satisfies both the Lipschitz and linear growth conditions across the intervals $[0, a_{1}]$ and $[ a_{1}, a_{2}]$. In the section $[a_{2}, a]$, we employ an alternative approach.\\
The situation for the interval $[ a_{2}, a]$ is covered first. For a given starting population size $U(t) \in R_{+}^2$, we observed that the parameters mentioned in the model are positive constants, implying their Lipschitz continuity.
 Therefore, it implies the existence of a unique local solution $U(t)\in [ a_{2}, \delta_{e}]$, where $\delta_{e}$ is defined as explosion-time. The obtained local solution is global if we establish that $\delta_{e}= \infty$. For this, we let $m_{0} \in R_{+}$ such that $U(a_{2}) \in \left[\frac{1}{m_{0}},m_{0}\right]$. Further, for each $m \ge m_{0}$, we define
	\begin{align*}
		\delta_{m}=\bigg\{t\in [a_{2}, \delta_{e}]:\frac{1}{m}\ge min\{U(t)\} \textbf{~~or}~~ max\{U(t)\} \ge m	\bigg\}
	\end{align*}
	One may easily observe that $ \delta_{m}$ defines a monotonically increasing sequence which will be bounded above. Therefore, it converges. Let $\lim_{m \to \infty} \delta_{m}=\delta_{\infty}$. Certainly $ \delta_{e} \ge \delta_{\infty}$. To prove that $\delta_{e}=\infty$ and that $U(t) \in R_{+}^2$ is the global solution, it is sufficient to establish that $\delta_{\infty}=\infty$, and that will imply $\delta_{e}=\infty$ for each $t \ge 0$. On the contrary, we assume that $\delta_{\infty}\neq \infty$. Then one may find  $a>0$ and $\zeta \in (0,1)$ such that \[P\{a \ge \delta_{\infty}\}> \zeta.\] 
	Next, we introduce a function $\tilde{\lambda}:R_{+}^2 \to R_{+}$ in $\lambda \in C^2$ space such that \[\tilde \lambda(U(t))=\sum_{i=1}^{2}X_{i}-2- \sum_{i=1}^{2} logX_{i}\]
	Since $z-1-logz \ge 0$ for each $z>0$, therefore, $\tilde \lambda \ge0$. Moreover, we assume that $m>m_{0}$ and $a>0$, therefore Ito formula implies the following
	\begin{align}\label{ito1}
		d \tilde \lambda(U)&=\mathbb{H}(U)+\sigma_{j}(X_{j}-1)d \mathbb{B}_{j}(t),
	\end{align}
	where
	\begin{align}
		\mathbb{H}(U) =\sum_{j=1}^{2}\left(1-\frac{1}{X_{j}}\right)X_{j}^{'}=\left(1-\frac{1}{r}\right)\bigg(-\rho_{1}r+\omega_{1}s(1-\epsilon s^2)+\psi_{1}\bigg)+\left(1-\frac{1}{s}\right)\bigg(-\rho_{2}s+\omega_{2}r(1-\epsilon r^2)+\psi_{2}\bigg)+ \sum_{j=1}^{2} \frac{\sigma_{j}^2}{2} \nonumber\\ = \psi_1+\rho_{1}+\psi_2+\rho_2-\rho_{1}s+w_1\bigg(s-\frac{s}{r}\bigg)- \epsilon\bigg(w_1^3+w_2r^3+w_1\frac{s^3}{r}+w_2\frac{r^3}{s} \bigg)-\frac{\psi_1}{r}-\rho_2s+w_2\bigg(r-\frac{r}{s}\bigg)-\frac{\psi_2}{s}+\sum_{j=1}^{2} \frac{\sigma_{j}^2}{2}\nonumber
	\end{align}
	Then, one may verify that \[\mathbb{H}(U)<\psi_1+\rho_{1}+\psi_2+\rho_2=\tilde \phi.\] Integrating \eqref{ito1} from $0$ to $\delta_{m} \wedge a$ and using the above fact, we deduce that  
\begin{align}
	E [\tilde \lambda(\delta_{m} \wedge a)]\le \tilde \lambda (U(a_{2}))+E\left[\int_{0}^{\delta_{m} \wedge a} \tilde \phi\right] \le  \tilde \lambda (U(a_{2}))+\tilde \phi a
\end{align}
We denote $\Pi_{m}=\{a>\delta_{m}\}$ for $m>m_0$. It is easy to see $P(\Pi_{m}) \ge \zeta$. Further, for every $\Omega \in \Pi_{m}$, one can find $U(\delta_{m},\omega)$ such that either  $U=\frac{1}{m}$ or $U=m$. These facts imply the following result\\
\[\left(logm+\frac{1}{m}-1\right) \wedge E(m-logm-1)<\tilde \lambda (U(\delta_{m})).\]
Therefore, we have
\begin{align}
	\tilde \lambda (U(a_{2}))+\tilde \phi a> E(1_{\Pi_{m}} \lambda (U(\delta_{m}))) \ge \sigma((m-logm-1) \wedge (logm+\frac{1}{m}-1) ),
\end{align}
where $1_{\Pi_{m}}$ denotes an indicator function of $\Pi$. On letting $m \to \infty$, we deduce that \[\infty>\tilde \lambda(U(a_{2}))+\tilde \phi a=\infty.\] The assumption is contradicted by the aforementioned fact. This means that $\delta_{\infty}=\infty$, which completes the evidence supporting the existence of a global solution.\\
The situation involving $[0,a_{2}]$ is then covered. We will move forward in this manner: Regarding each $t \in [0,a_2]$, we get
\begin{align}
	|e_{1}(r,s,t)|^2 &\le |-\rho_{1}r+\omega_{1}s(1-\epsilon s^2)+\psi_{1}|^2 
	\le 4[\rho_{1}^2|r|^2+\omega_{1}^2|s|^2+\epsilon^2\omega_{1}^2|s^2|^2+ \psi_{1}^2]\nonumber \\
	&\le 4[\rho_{1}^2|r|^2+\omega_{1}^2  \sup_{t\in [0,T]} |s|^2+\epsilon^2\omega_{1}^2  \sup_{t\in [0,T]} |s^2|^2+ \psi_{1}^2]
	\le 4\bigg[\rho_{1}^2|r|^2+\omega_{1}^2||s^2||_{\infty}+\epsilon^2\omega_{1}^2||s^4||_{\infty}+ \psi_{1}^2\bigg] \nonumber\\
	&\le 4 \bigg(\omega_{1}^2||s^2||_{\infty}+\epsilon^2\omega_{1}^2||s^4||_{\infty}+ \psi_{1}^2\bigg) \left[1+\frac{\rho_{1}^2|r|^2}{\bigg(\omega_{1}^2||s^2||_{\infty}+\epsilon^2\omega_{1}^2||s^4||_{\infty}+ \psi_{1}^2\bigg)}\right] \le 4M_{1}\left(1+\frac{\rho_{1}^2}{M_1} |r|^2\right)\le c_{1}(1+|r|^2)
\end{align} 
with $c_1=4M_1$ and $ M_{1}=\omega_{1}^2||s^2||_{\infty}+\epsilon^2\omega_{1}^2||s^4||_{\infty}+ \psi_{1}^2$ provided $\frac{\rho_{1}^2}{M_1}<1$. 

\begin{align}
	|e_{2}(r,s,t)|^2 &\le |-\rho_{2}s+\omega_{2}r(1-\epsilon r^2)+\psi_{2}|^2 
	\le 4[\rho_{2}^2|s|^2+\omega_{2}^2|r|^2+\epsilon^2\omega_{2}^2|r^2|^2+ \psi_{2}^2]\nonumber \\
	&\le 4[\rho_{2}^2|s|^2+\omega_{2}^2  \sup_{t\in [0,a]} |r|^2+\epsilon^2\omega_{2}^2  \sup_{t\in [0,a]} |r^2|^2+ \psi_{2}^2]
	\le 4\bigg[\rho_{2}^2|s|^2+\omega_{2}^2||r^2||_{\infty}+\epsilon^2\omega_{2}^2||r^4||_{\infty}+ \psi_{2}^2\bigg] \nonumber\\
	&\le 4 \bigg(\omega_{2}^2||r^2||_{\infty}+\epsilon^2\omega_{2}^2||r^4||_{\infty}+ \psi_{2}^2\bigg) \left[1+\frac{\rho_{2}^2|s|^2}{\bigg(\omega_{2}^2||r^2||_{\infty}+\epsilon^2\omega_{2}^2||r^4||_{\infty}+ \psi_{2}^2\bigg)}\right] \le 4M_{1}\left(1+\frac{\rho_{2}^2}{M_2} |s|^2\right)\le c_{2}(1+|s|^2)
\end{align} 
with $c_2=4M_2$ and $ M_{2}=\omega_{2}^2||r^2||_{\infty}+\epsilon^2\omega_{2}^2||r^4||_{\infty}+ \psi_{2}^2$ provided $\frac{\rho_{2}^2}{M_2}<1$.\\
It may be inferred from the foregoing facts that $e_1$ and $e_2$ meet the linear growth criteria on $[0,a_2]$. Additionally, the following conclusions can be drawn for every $t\in [0,a_2]$.
\begin{align}
	|e_{1}(t,r_{1},s)-e_{1}(t,r_{2},s)|^2 \le \rho_{1}^2|r_{1}-r_{2}|^2 \le \tilde m_{1}|r_{1}-r_{2}|^2; \tilde m_{1}=\rho_1^2\nonumber \\
	|e_{2}(t,r,s_1)-e_{2}(t,r,s_2)|^2 \le \rho_{2}^2|s_{1}-s_{2}|^2 \le \tilde m_{2}|s_{1}-s_{2}|^2; \tilde m_{2}=\rho_2^2. \nonumber
\end{align}
Regarding $[0,a]$, if we assume that the solutions are positive, then $\max\left\{\frac{\rho_{1}^2}{M_1},\frac{\rho_{2}^2}{M_2}\right\}<1$ is the set of values for which the piecewise system has a positive solution. A similar result for the second piecewise love dynamical model's existence can be demonstrated.

	\section{Numerical algorithm}
Numerical solutions are established for the proposed systems in this section. These are obtained by the Newton interpolation, and applied in first and second models $[0,H]$ is distributed as following:
\begin{align*}
	0\leq t_{0} \leq t_{1} \leq t_{2}...\leq t_{k_{1}} =\\ a_{1}\leq t_{k_{1}+1} \leq t_{k_{1}+2} \leq t_{k_{1}+3}...\leq t_{k_{2}}=a.
\end{align*}
$h(t,\xi)$ will be approximated by Newton interpolation \label{Newton_method} $P(t)$, which will be represented as:
\begin{align}\label{Newton_method}
	e(t,U) \approx P(t)=e(t_{j-1},U_{j-1})+ \frac{e(t_{j-1},U_{j-1})-e(t_{j-2},U_{j-2})}{\Delta t}(t-t_{j-2})+\\ \frac{e(t_{j},U_{j})-2e(t_{j-1},U_{j-1})+e(t_{j-2},U_{j-2})}{2(\Delta t)^2}(t-t_{j-1})(t-t_{j-2}).  \nonumber
\end{align} 
\subsection{Numerical algorithm in the aspect of Caputo derivative}	
Numerical solutions for the model applying classical power law can be provided as
\begin{align}
	\begin{cases}
		&U^{k_{1}} = U(0)+\frac{1}{12}\sum_{j_1 = 2}^{k_{1}}\bigg[23e(t_{j_1}, U(t_{j_1}))-16e(t_{j_1 - 1}, U(t_{j_1 - 1}))+5e(t_{j_1 - 2}, U(t_{j_1 - 2})) \bigg]* \Delta t, 0 \le t\le a_{1}\\
		&U^{k_{2}} = U(a_{1})
		+ \frac{h^\alpha}{\Gamma{(\alpha)}} \sum_{j_2 = k_{1}+3}^{k_{2}}\bigg(e(t_{j_2 - 1}, U(t_{j_2 -1}))\bigg[\frac{(k_2-j_2+1)^{\alpha}-(k_2-j_2)^{\alpha} }{\alpha} \bigg]\\
		&+(e(t_{j_2 - 1}, U(t_{j_2 - 1}))-e(t_{j_2 - 2}, U(t_{j_2 - 2}))) \bigg[(k_2-j_2+2)\frac{(k_2-j_2+1)^{\alpha}-(k_2-j_2)^{\alpha}}{\alpha} 
		-\frac{(k_2-j_2+1)^{\alpha+1}-(k_2-j_2)^{\alpha+1} }{\alpha+1}\bigg]\\
		&+\frac{e(t_{j_2}, U(t_{j_2}))-2e(t_{j_2-1 }, U(t_{j_2-1 }))+e(t_{j_2 -2}, U(t_{j_2 -2}))}{2} \bigg[(k_2-j_2+2)(k_2-j_2+1)\frac{(k_2-j_2+1)^{\alpha}-(k_2-j_2)^{\alpha} }{\alpha}\\
		&-(2k_2-2j_2+3)\frac{(k_2-j_2+1)^{\alpha+1}-(k_2-j_2)^{\alpha+1} }{\alpha+1}
		+\frac{(k_2-j_2+1)^{\alpha+2}-(k_2-j_2
			)^{\alpha+2} }{\alpha+2}\bigg]\bigg),
	~a_{1} \le t\le a_2 \\
		&U^{k_{3}}= U(a_{2})+ \frac{1}{12}\sum_{j_3 = k_{2}+3}^{k_{3}} \bigg[23e(t_{j_3}, U(t_{j_3}))-16e(t_{j_3 - 1}, U(t_{j_3 - 1})) +5e(t_{j_3 - 2}, U(t_{j_3 - 2}))\bigg]*\Delta t+ \\
		&~\sigma_{i}\sum_{j_3 = k_{2}+3}^{k_{3}} U(\mathbb{B}_{i}(t_{j_3})-\mathbb{B}_{i}(t_{j_3-1})), a_{2} \le t\le a. 
	\end{cases}
\end{align}
\subsection*{Computer simulation of first and second modified love dynamical systems in the aspect of Caputo derivative:}  
\begin{figure}[H]
	\centering
	\begin{subfigure}{0.325\textwidth}
		\includegraphics[width=\textwidth]{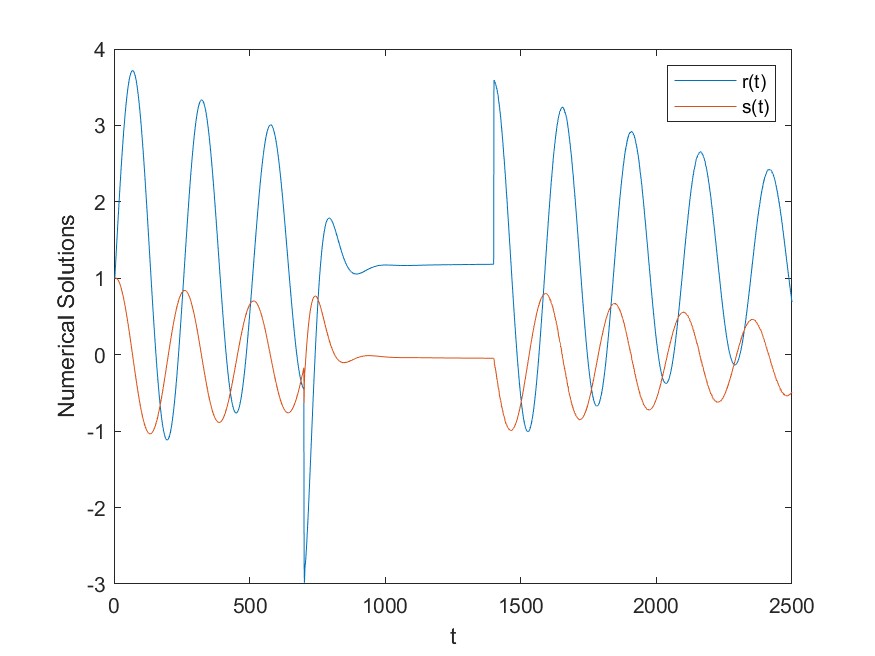}
		\caption{$\alpha=0.79$}
	\end{subfigure}
	\hfill
	\begin{subfigure}{0.325\textwidth}
		\includegraphics[width=\textwidth]{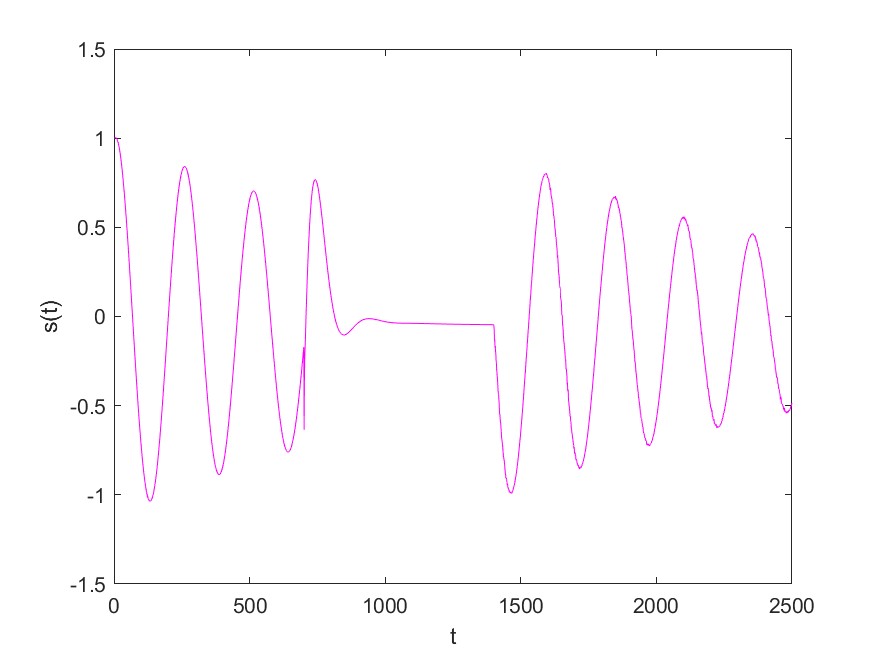}
		\caption{$\alpha=0.79$}
	\end{subfigure}
	\hfill
	\begin{subfigure}{0.325\textwidth}
		\includegraphics[width=\textwidth]{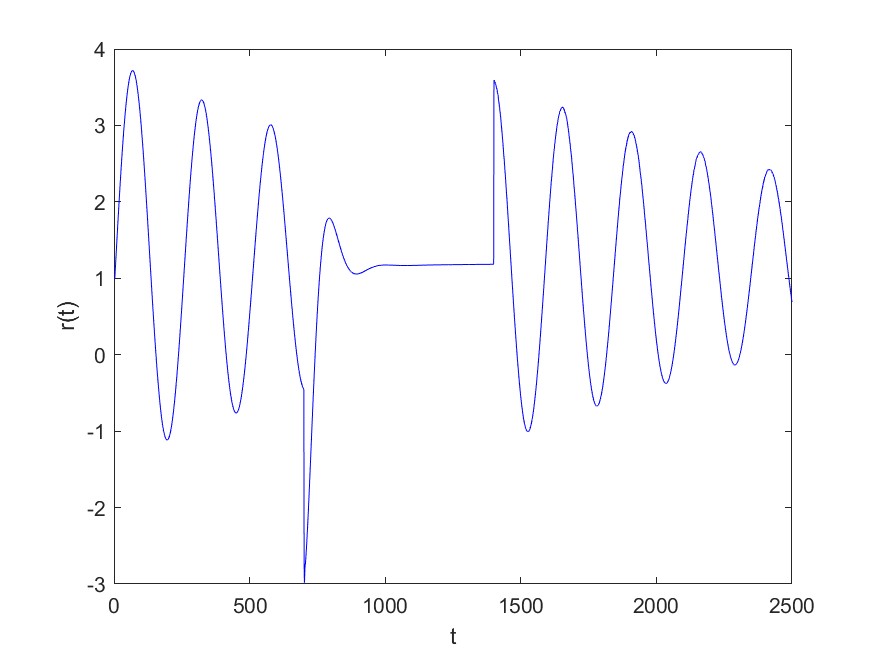}
		\caption{$\alpha=0.79$}
	\end{subfigure}
	\hfill
	\begin{subfigure}{0.325\textwidth}
		\includegraphics[width=\textwidth]{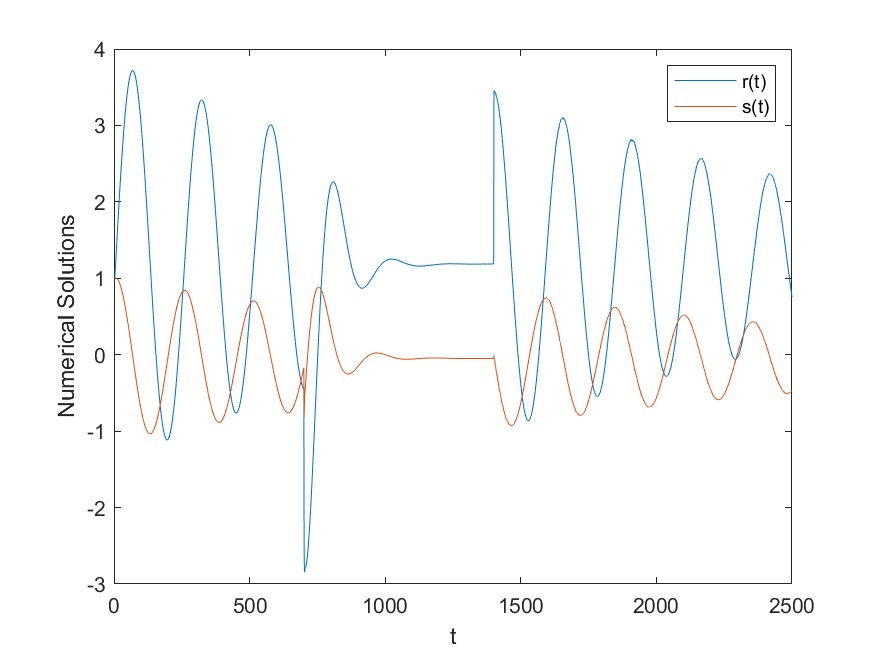}
		\caption{$\alpha=0.85$}
	\end{subfigure}
	\hfill
	\begin{subfigure}{0.325\textwidth}		\includegraphics[width=\textwidth]{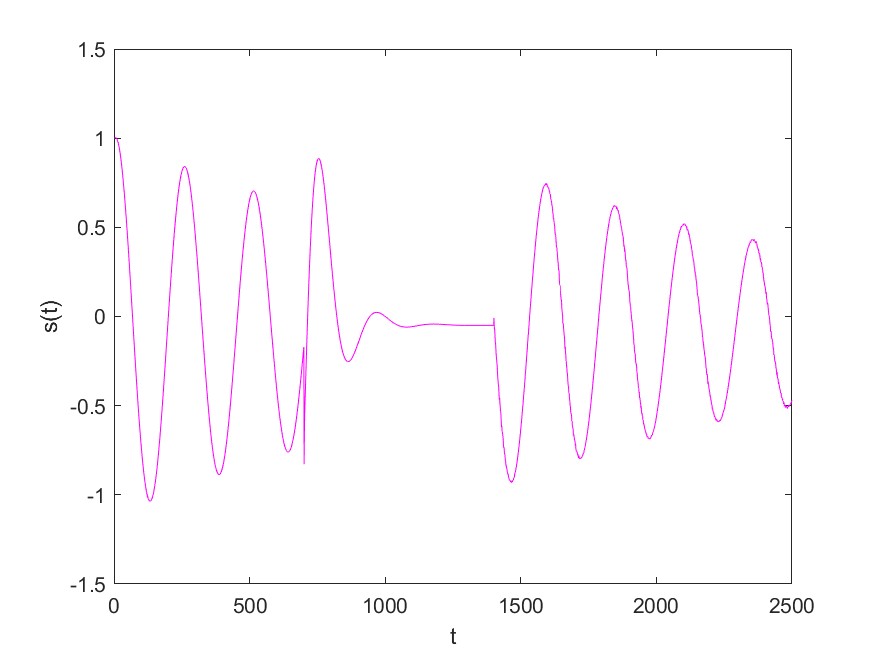}
		\caption{$\alpha=0.85$}
	\end{subfigure}
	\hfill
	\begin{subfigure}{0.325\textwidth}
		\includegraphics[width=\textwidth]{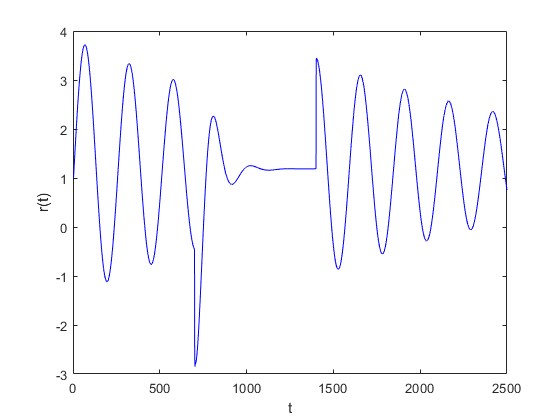}
		\caption{$\alpha=0.85$}
	\end{subfigure}
	\hfill
	\begin{subfigure}{0.325\textwidth}
		\includegraphics[width=\textwidth]{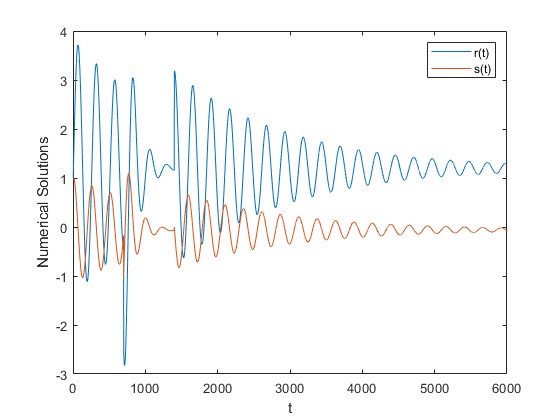}
		\caption{$\alpha=0.92$.}
	\end{subfigure}
	\hfill
	\begin{subfigure}{0.325\textwidth}
		\includegraphics[width=\textwidth]{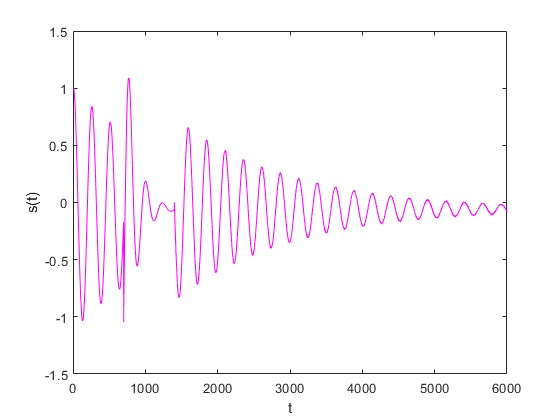}
		\caption{$\alpha=0.92$.}
	\end{subfigure}
	\hfill
	\begin{subfigure}{0.325\textwidth}
		\includegraphics[width=\textwidth]{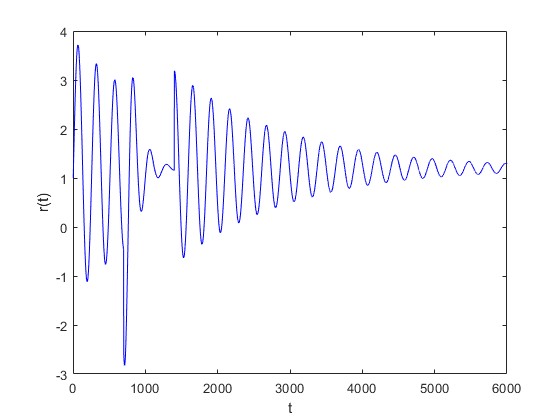}
		\caption{$\alpha=0.92$.}
	\end{subfigure}
	\hfill
	\begin{subfigure}{0.325\textwidth}
		\includegraphics[width=\textwidth]{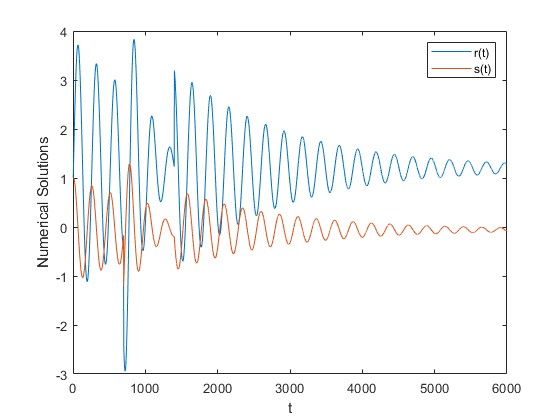}
		\caption{$\alpha=0.97$.}
	\end{subfigure}
	\hfill
	\begin{subfigure}{0.325\textwidth}
		\includegraphics[width=\textwidth]{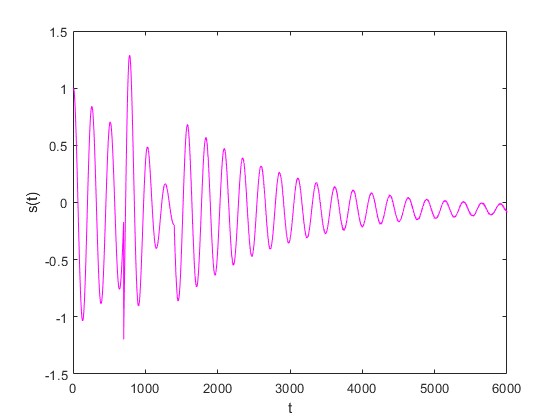}
		\caption{$\alpha=0.97$.}
	\end{subfigure}
	\hfill
	\begin{subfigure}{0.325\textwidth}
		\includegraphics[width=\textwidth]{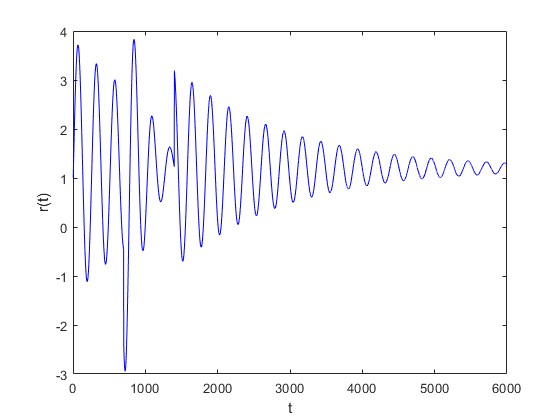}
		\caption{$\alpha=0.97$.}
	\end{subfigure}

	\caption{The Numerical simulations for first model \eqref{caputolinear} with  $t(0)=0$, $s(0)=1$, $p(0)=1$, $\rho_{1}=0.12$, $\rho_{2}=0.05$, $\psi_{1}=0.8$, $\psi_{2}=0.81$, $\gamma_{1}=0.5$, $\gamma_{2}=1.2$, $\sigma_{1}=0.02$, $\sigma_{2}=0.01$, $\omega_{1}=6.1$, $\omega_{2}=-1$.}
	\label{fig:fig1l}
\end{figure}

\begin{figure}[H]
	\centering
	\begin{subfigure}{0.32\textwidth}
		\includegraphics[width=\textwidth]{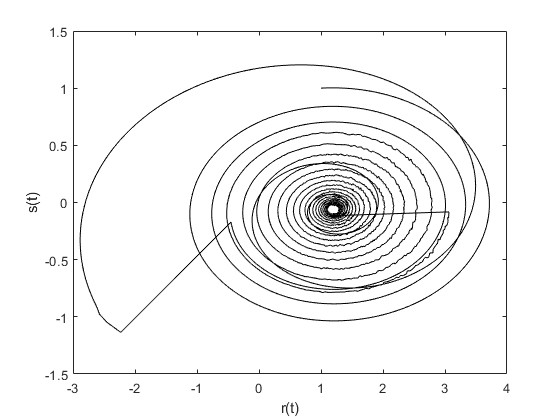}
		\caption{$\rho_{2}=0.01$, $\omega_{2}=-1$.}
	\end{subfigure}
	\hfill
	\begin{subfigure}{0.32\textwidth}
		\includegraphics[width=\textwidth]{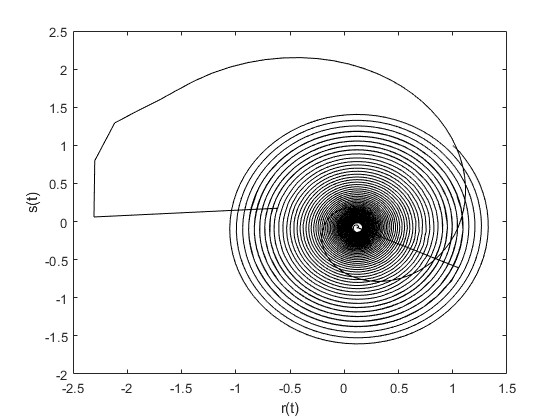}
		\caption{$\rho_{2}=0.01$, $\omega_{2}=-10$.}
	\end{subfigure}
	\hfill
	\begin{subfigure}{0.32\textwidth}
		\includegraphics[width=\textwidth]{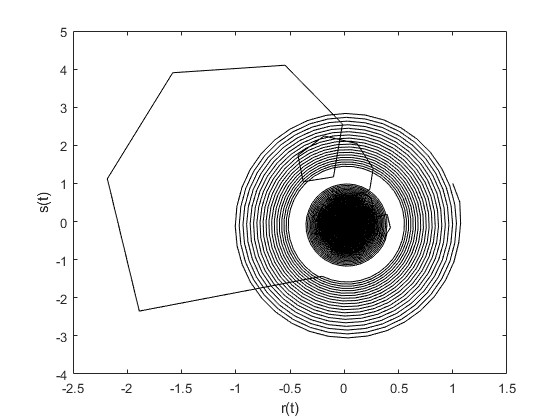}
		\caption{$\rho_{2}=0.01$, $\omega_{2}=-50$.}
	\end{subfigure}
	\hfill
	\begin{subfigure}{0.32\textwidth}
		\includegraphics[width=\textwidth]{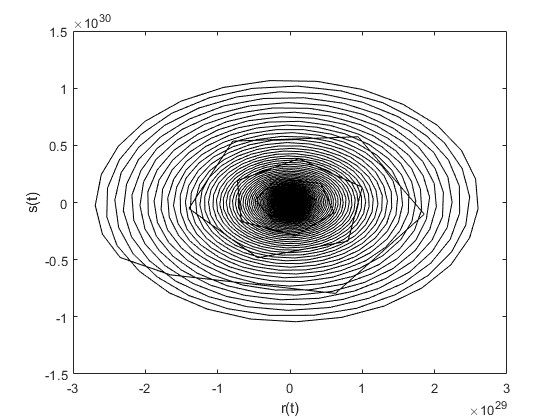}
		\caption{$\rho_{2}=0.01$, $\omega_{2}=-100$.}
	\end{subfigure}
	\hfill
	\begin{subfigure}{0.32\textwidth}
		\includegraphics[width=\textwidth]{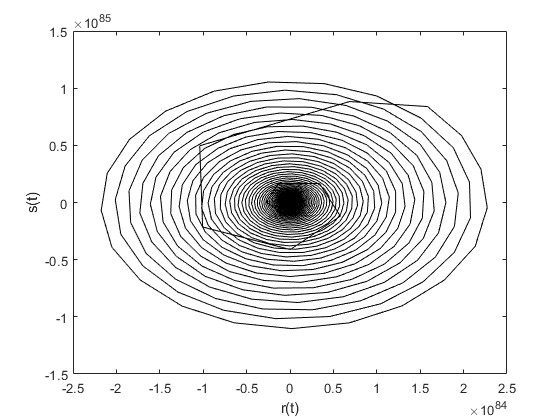}
		\caption{$\rho_{2}=0.01$, $\omega_{2}=-150$.}
	\end{subfigure}
	\hfill
	\begin{subfigure}{0.32\textwidth}
		\includegraphics[width=\textwidth]{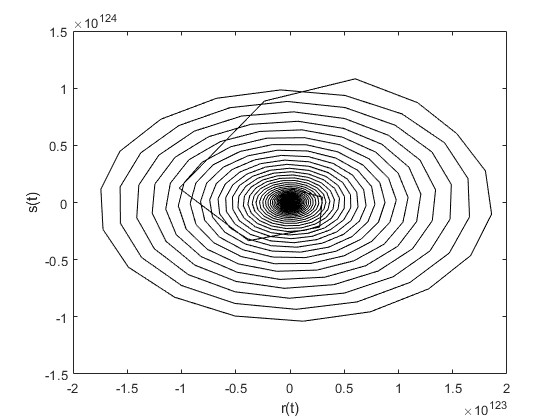}
		\caption{$\rho_{2}=0.01$, $\omega_{2}=-200$.}
	\end{subfigure}
	\hfill
	\begin{subfigure}{0.32\textwidth}
		\includegraphics[width=\textwidth]{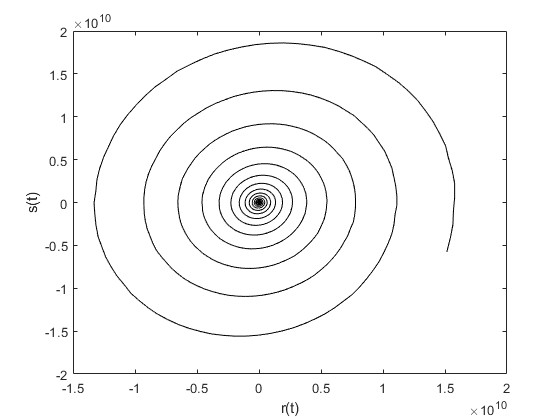}
		\caption{$\rho_{2}=-1$, $\omega_{2}=-10$.}
	\end{subfigure}
	\hfill
	\begin{subfigure}{0.32\textwidth}
		\includegraphics[width=\textwidth]{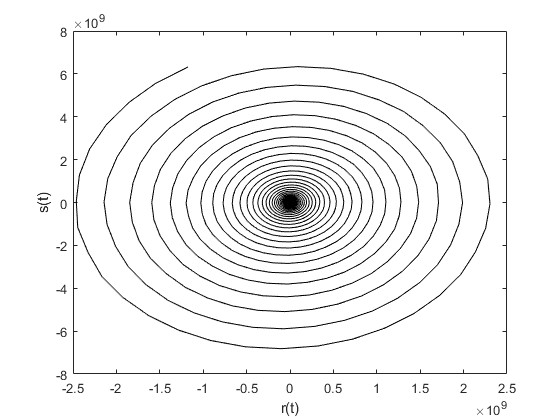}
		\caption{$\rho_{2}=-1$, $\omega_{2}=-50$.}
	\end{subfigure}
	\hfill
	\begin{subfigure}{0.32\textwidth}
		\includegraphics[width=\textwidth]{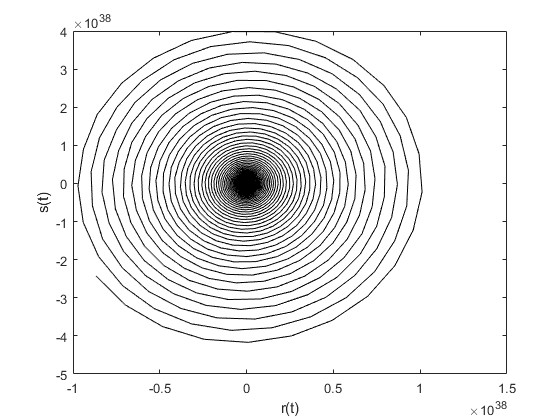}
		\caption{$\rho_{2}=-1$, $\omega_{2}=-100$.}
	\end{subfigure}
	\hfill
	\begin{subfigure}{0.32\textwidth}
		\includegraphics[width=\textwidth]{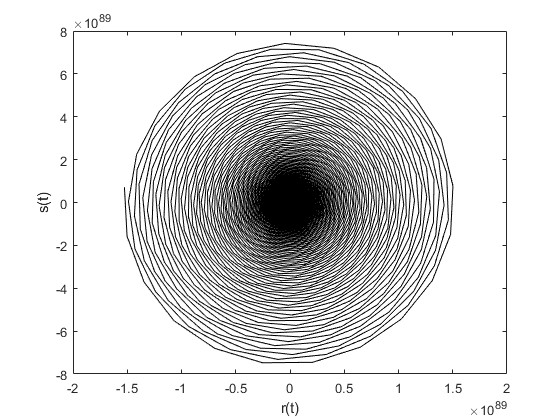}
		\caption{$\rho_{2}=-1$, $\omega_{2}=-150$.}
	\end{subfigure}
	\hfill
	\begin{subfigure}{0.32\textwidth}
		\includegraphics[width=\textwidth]{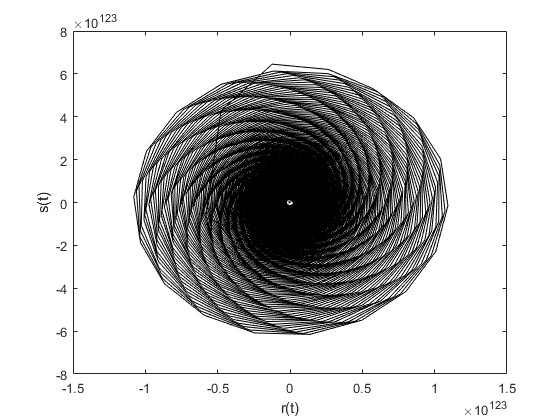}
		\caption{$\rho_{2}=-1$, $\omega_{2}=-200$.}
	\end{subfigure}
	\hfill
	\begin{subfigure}{0.32\textwidth}
		\includegraphics[width=\textwidth]{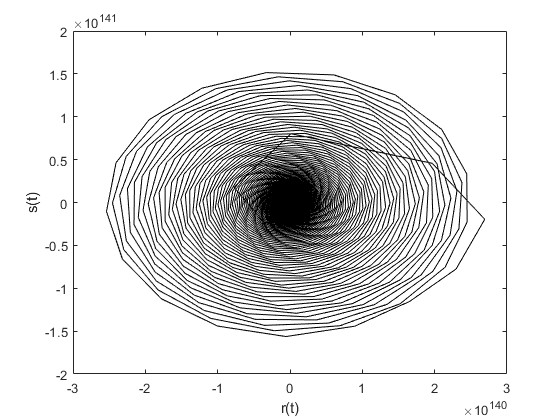}
		\caption{$\rho_{2}=-1$, $\omega_{2}=-225$.}
	\end{subfigure}
	\caption{Chaotic dynamics for the first model \eqref{caputolinear} with  $t(0)=0$, $s(0)=1$, $p(0)=1$, $\rho_{1}=0.12$, $\psi_{1}=0.8$, $\psi_{2}=0.81$, $\gamma_{1}=0.5$, $\gamma_{2}=1.2$, $\sigma_{1}=0.02$, $\sigma_{2}=0.01$, $\omega_{1}=6.1$, and $\alpha=0.95$..}
	\label{fig:fig2l}
\end{figure}

\begin{figure}[H]
	\centering
	\begin{subfigure}{0.325\textwidth}
		\includegraphics[width=\textwidth]{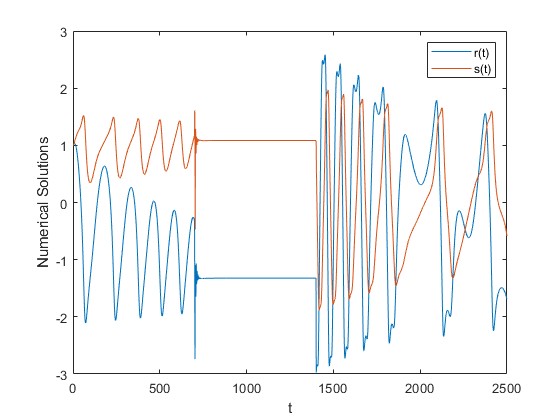}
		\caption{$\alpha=0.8$, $\epsilon=1$}
	\end{subfigure}
	\hfill
	\begin{subfigure}{0.325\textwidth}
		\includegraphics[width=\textwidth]{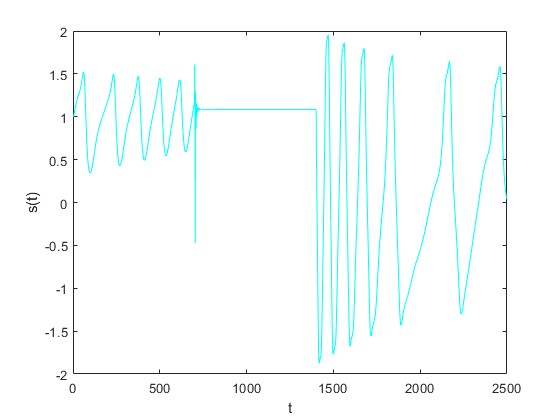}
		\caption{$\alpha=0.8$, $\epsilon=1$.}
	\end{subfigure}
	\hfill
	\begin{subfigure}{0.325\textwidth}
		\includegraphics[width=\textwidth]{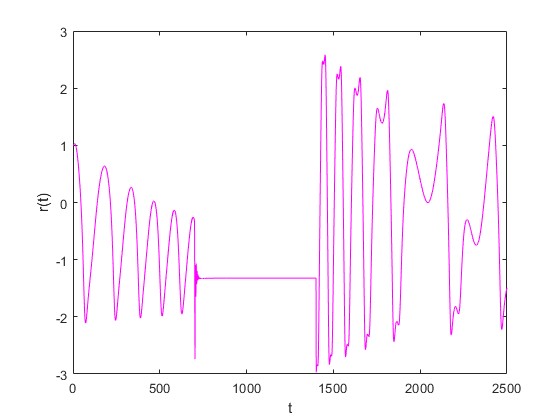}
		\caption{$\alpha=0.1$, $\epsilon=1$.}
	\end{subfigure}
	\hfill
	\begin{subfigure}{0.325\textwidth}
		\includegraphics[width=\textwidth]{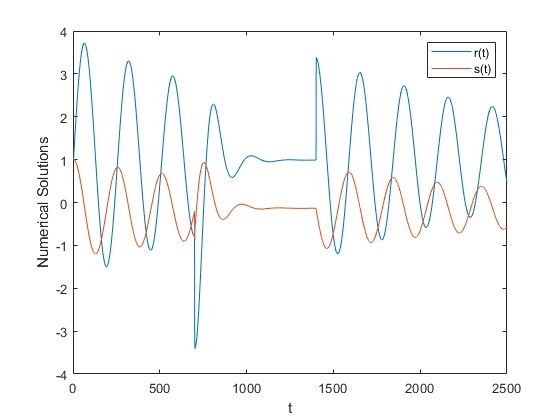}
		\caption{$\alpha=0.86$, $\epsilon=0$.}
	\end{subfigure}
	\hfill
	\begin{subfigure}{0.325\textwidth}		\includegraphics[width=\textwidth]{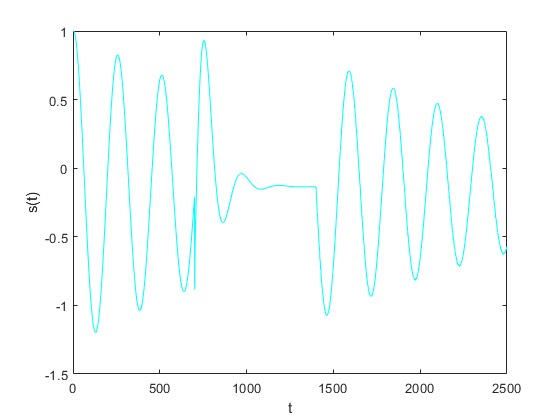}
		\caption{$\alpha=0.86$, $\epsilon=0$.}
	\end{subfigure}
	\hfill
	\begin{subfigure}{0.325\textwidth}
		\includegraphics[width=\textwidth]{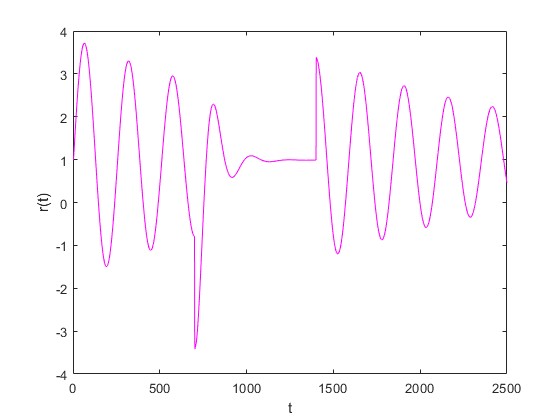}
		\caption{$\alpha=0.86$, $\epsilon=0$.}
	\end{subfigure}
	\hfill
	\begin{subfigure}{0.325\textwidth}
		\includegraphics[width=\textwidth]{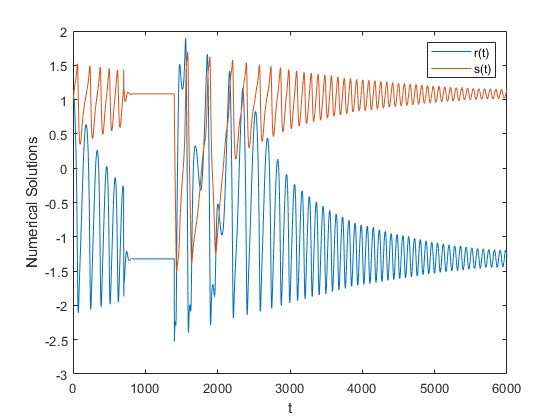}
		\caption{$\alpha=0.93$, $\epsilon=1$.}
	\end{subfigure}
	\hfill
	\begin{subfigure}{0.325\textwidth}
		\includegraphics[width=\textwidth]{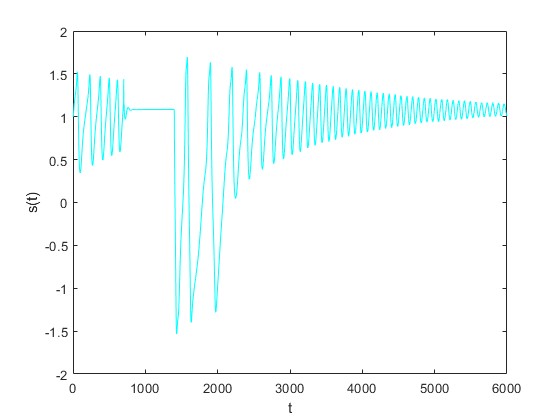}
		\caption{$\alpha=0.93$, $\epsilon=1$.}
	\end{subfigure}
	\hfill
	\begin{subfigure}{0.325\textwidth}
		\includegraphics[width=\textwidth]{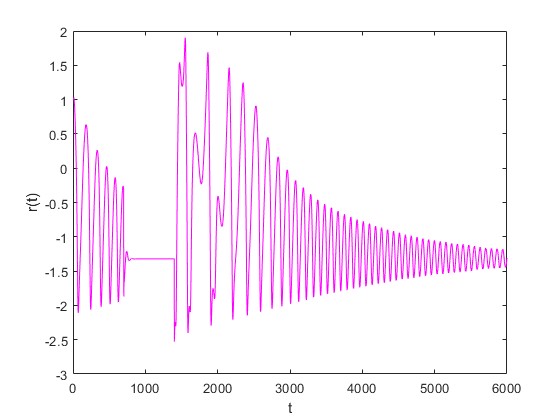}
		\caption{$\alpha=0.93$, $\epsilon=1$.}
	\end{subfigure}
	\hfill
	\begin{subfigure}{0.325\textwidth}
		\includegraphics[width=\textwidth]{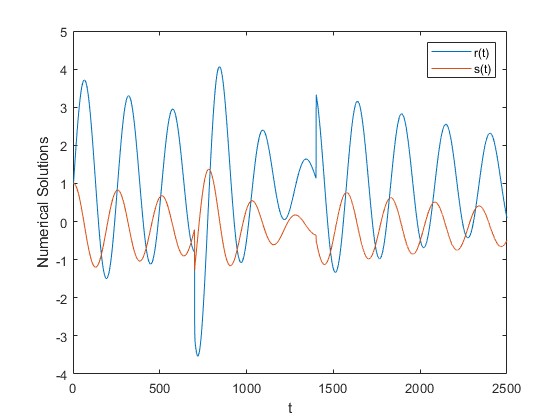}
		\caption{$\alpha=0.98$, $\epsilon=0$.}
	\end{subfigure}
	\hfill
	\begin{subfigure}{0.325\textwidth}
		\includegraphics[width=\textwidth]{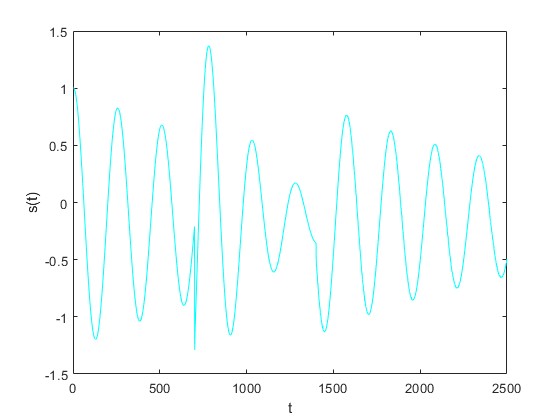}
		\caption{$\alpha=0.98$, $\epsilon=0$.}
	\end{subfigure}
	\hfill
	\begin{subfigure}{0.325\textwidth}
		\includegraphics[width=\textwidth]{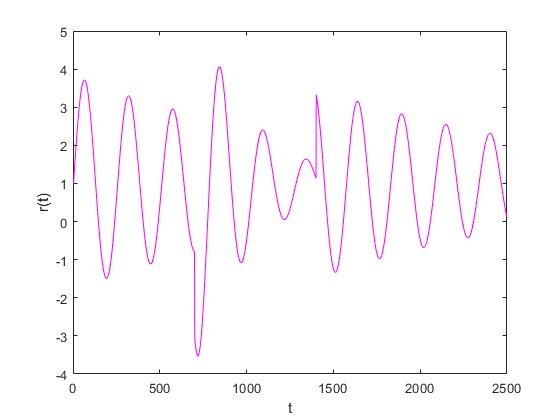}
		\caption{$\alpha=0.98$, $\epsilon=0$.}
	\end{subfigure}
	\caption{The Numerical simulations for second model \eqref{caputononlinear} with $t(0)=0$, $r(0)=1$, $s(0)=1$,  $\rho_{1}=0.12$, $\rho_{2}=0.01$,
		$\phi_{1}=1$, $\phi_{2}=1$,
		$\sigma_{1}=0.01$, $\sigma_{2}=0.02$,
		$\omega_{1}=6.1$, $\omega_{2}=-1$.}
	\label{fig:fig3l}
\end{figure}

\begin{figure}[H]
	\centering
	\begin{subfigure}{0.32\textwidth}
		\includegraphics[width=\textwidth]{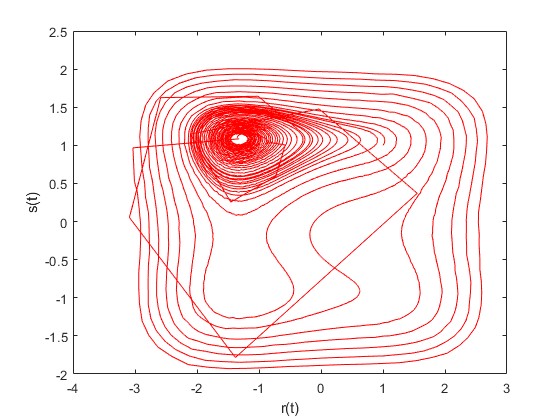}
		\caption{$\alpha=0.78$, $\epsilon=1$.}
	\end{subfigure}
	\hfill
	\begin{subfigure}{0.32\textwidth}
		\includegraphics[width=\textwidth]{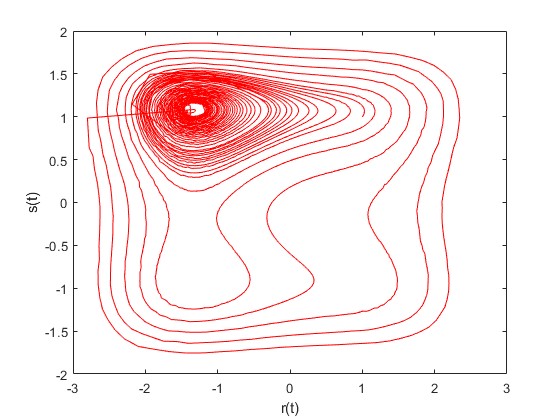}
		\caption{$\alpha=0.85$, $\epsilon=1$.}
	\end{subfigure}
	\hfill
	\begin{subfigure}{0.32\textwidth}
		\includegraphics[width=\textwidth]{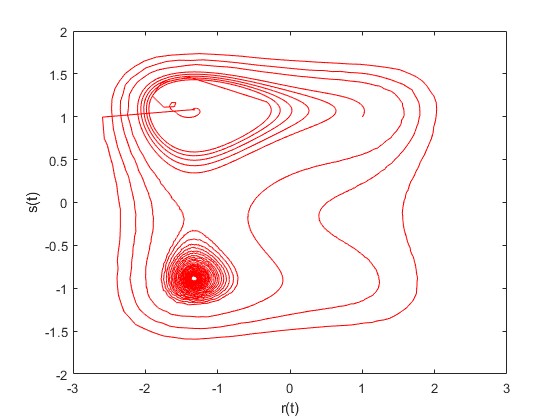}
		\caption{$\alpha=0.91$, $\epsilon=1$.}
	\end{subfigure}
	\hfill
	\begin{subfigure}{0.32\textwidth}
		\includegraphics[width=\textwidth]{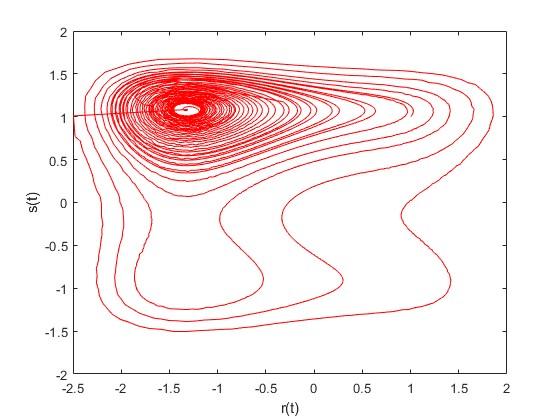}
		\caption{$\alpha=0.94$, $\epsilon=1$.}
	\end{subfigure}
	\hfill
	\begin{subfigure}{0.32\textwidth}
		\includegraphics[width=\textwidth]{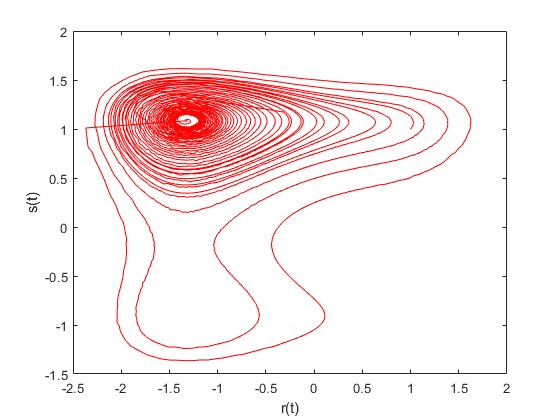}
		\caption{$\alpha=0.98$, $\epsilon=1$.}
	\end{subfigure}
	\hfill
	\begin{subfigure}{0.32\textwidth}
		\includegraphics[width=\textwidth]{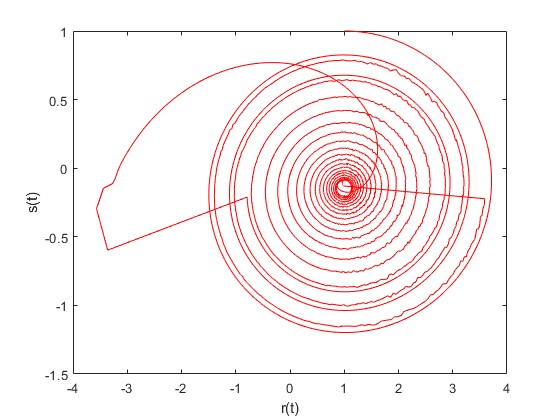}
		\caption{$\alpha=0.78$, $\epsilon=0$.}
	\end{subfigure}
	\hfill
	\begin{subfigure}{0.32\textwidth}
		\includegraphics[width=\textwidth]{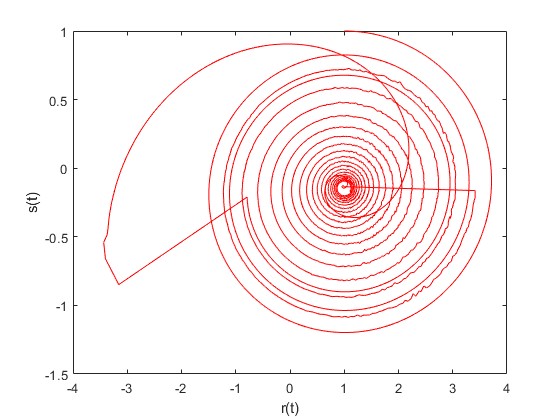}
		\caption{$\alpha=0.85$, $\epsilon=0$.}
	\end{subfigure}
	\hfill
	\begin{subfigure}{0.32\textwidth}
		\includegraphics[width=\textwidth]{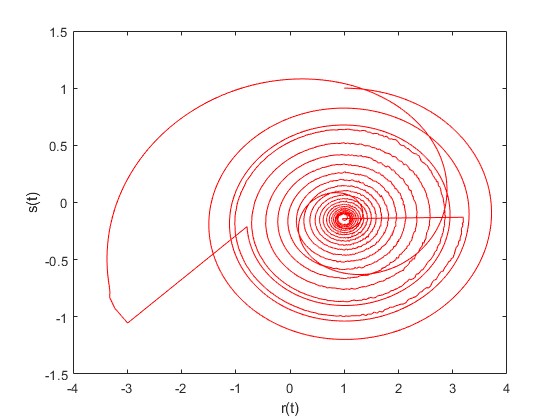}
		\caption{$\alpha=0.91$, $\epsilon=0$.}
	\end{subfigure}
	\hfill
	\begin{subfigure}{0.32\textwidth}
		\includegraphics[width=\textwidth]{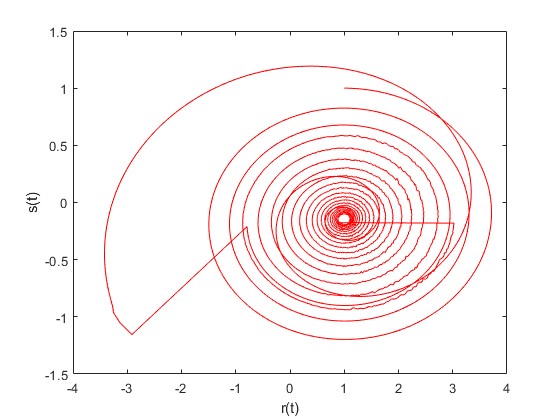}
		\caption{$\alpha=0.94$, $\epsilon=0$.}
	\end{subfigure}
	\hfill
	\begin{subfigure}{0.32\textwidth}
		\includegraphics[width=\textwidth]{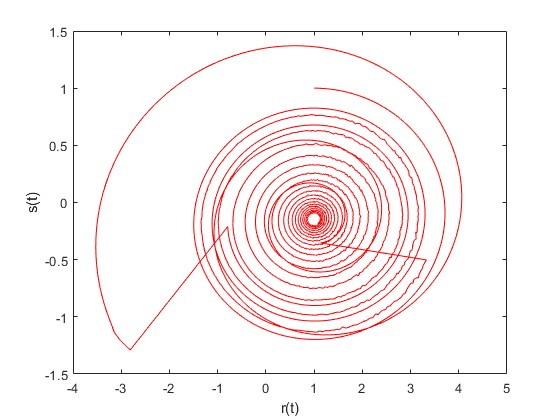}
		\caption{$\alpha=0.98$, $\epsilon=0$.}
	\end{subfigure}
	\hfill
	\begin{subfigure}{0.32\textwidth}
		\includegraphics[width=\textwidth]{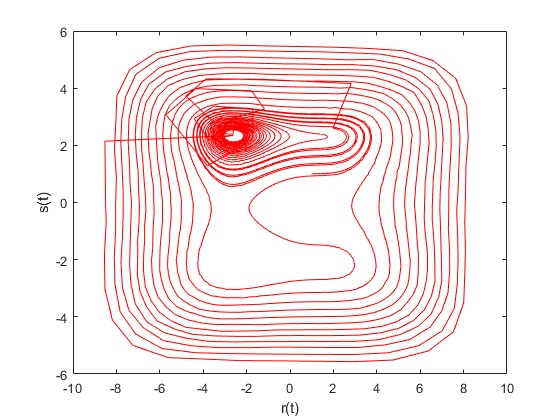}
		\caption{$\alpha=0.85$, $\epsilon=0.2$.}
	\end{subfigure}
	\hfill
	\begin{subfigure}{0.32\textwidth}
		\includegraphics[width=\textwidth]{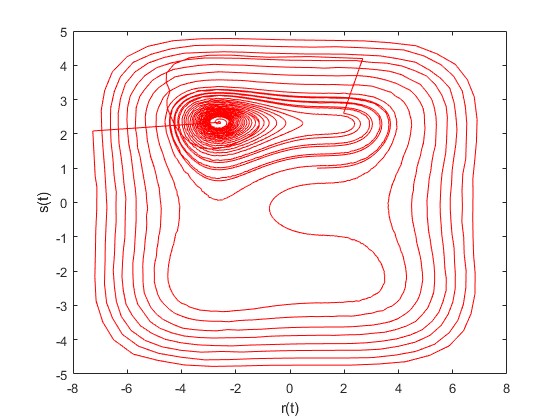}
		\caption{$\alpha=0.98$, $\epsilon=0.2$.}
	\end{subfigure}
	\caption{Chaotic dynamics for the second model \eqref{caputononlinear} with  $t(0)=0$, $r(0)=1$, $s(0)=1$,  $\rho_{1}=0.12$, $\rho_{2}=0.01$,
		$\phi_{1}=1$, $\phi_{2}=1$,
		$\sigma_{1}=0.01$, $\sigma_{2}=0.02$,
		$\omega_{1}=6.1$, $\omega_{2}=-1$.}
	\label{fig:fig4l}
\end{figure}

	It presents the numerical simulations and chaotic phase portrait via first \eqref{caputolinear}and second \eqref{caputononlinear} models for love dynamical system. The 
	figures \ref{fig:fig1l}(a)-\ref{fig:fig1l}(l), and \ref{fig:fig2l}(a)-\ref{fig:fig2l}(l) display the numerical simulations, while the figures \ref{fig:fig3l}(a)-\ref{fig:fig3l}(l), and \ref{fig:fig4l}(a)-\ref{fig:fig4l}(l) represent the chaotic behaviors of phase portrait for the first and second models in the aspect of Caputo derivative.

	\subsection{Numerical algorithm in the aspect of Atangana-Baleanu derivative:}
Numerical solutions for the model applying Mittag-Leffler law, can be provided as
\begin{align}
	\begin{cases}
		&U^{k_{1}} = U(0)+\frac{1}{12}\sum_{j_1 = 2}^{k_{1}}\bigg[23e(t_{j_1}, U(t_{j_1}))-16e(t_{j_1 - 1}, U(t_{j_1 - 1}))+5e(t_{j_1 - 2}, U(t_{j_1 - 2})) \bigg]* \Delta t, 0 \le t\le a_{1},  \\
		&U^{k_{2}} = U(a_{1})+
		\frac{1-\alpha}{AB(\alpha)} e(t_{k_2}, U(t_{k_2}))+\\
		&+\frac{\alpha}{AB(\alpha) \Gamma(\alpha+1)} \frac{h^\alpha}{\Gamma{(\alpha)}} \sum_{j_2 = k_{1}+3}^{k_{2}}\bigg(e(t_{j_2 - 1}, U(t_{j_2 -1}))\bigg[\frac{(k_2-j_2+1)^{\alpha}-(k_2-j_2)^{\alpha} }{\alpha} \bigg]\\
		&+(e(t_{j_2 - 1}, U(t_{j_2 - 1}))-e(t_{j_2 - 2}, U(t_{j_2 - 2}))) \bigg[(k_2-j_2+2)\frac{(k_2-j_2+1)^{\alpha}-(k_2-j_2)^{\alpha}}{\alpha} 
		-\frac{(k_2-j_2+1)^{\alpha+1}-(k_2-j_2)^{\alpha+1} }{\alpha+1}\bigg]\\
		&+\frac{e(t_{j_2}, U(t_{j_2}))-2e(t_{j_2-1 }, U(t_{j_2-1 }))+e(t_{j_2 -2}, U(t_{j_2 -2}))}{2} \bigg[(k_2-j_2+2)(k_2-j_2+1)\frac{(k_2-j_2+1)^{\alpha}-(k_2-j_2)^{\alpha} }{\alpha}\\
		&-(2k_2-2j_2+3)\frac{(k_2-j_2+1)^{\alpha+1}-(k_2-j_2)^{\alpha+1} }{\alpha+1}
		+\frac{(k_2-j_2+1)^{\alpha+2}-(k_2-j_2)^{\alpha+2} }{\alpha+2}\bigg]\bigg), 
	~a_{1} \le t\le a_2 \\
		&U^{k_{3}}= U(a_{2})+ \frac{1}{12}\sum_{j_3 = k_{2}+3}^{k_{3}} \bigg[23e(t_{j_3}, U(t_{j_3}))-16e(t_{j_3 - 1}, U(t_{j_3 - 1})) +5e(t_{j_3 - 2}, U(t_{j_3 - 2}))\bigg]*\Delta t+ \\
		&~\sigma_{i}\sum_{j_3 = k_{2}+3}^{k_{3}} U(\mathbb{B}_{i}(t_{j_3})-\mathbb{B}_{i}(t_{j_3-1})), a_{2} \le t\le a. 
	\end{cases}
\end{align}
\subsection*{Computer simulation of first and second modified love dynamical systems in the aspect of Atangana-Baleanu derivative:}  
\begin{figure}[H]
	\centering
	\begin{subfigure}{0.325\textwidth}
		\includegraphics[width=\textwidth]{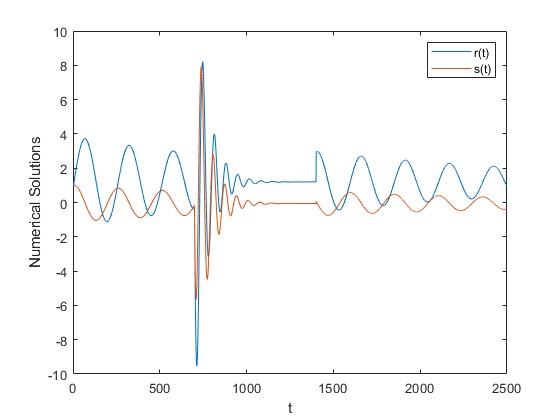}
		\caption{$\alpha=0.86$.}
	\end{subfigure}
	\hfill
	\begin{subfigure}{0.325\textwidth}
		\includegraphics[width=\textwidth]{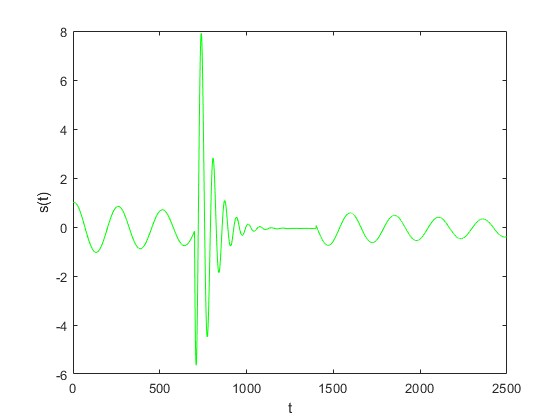}
		\caption{$\alpha=0.86$.}
	\end{subfigure}
	\hfill
	\begin{subfigure}{0.325\textwidth}
		\includegraphics[width=\textwidth]{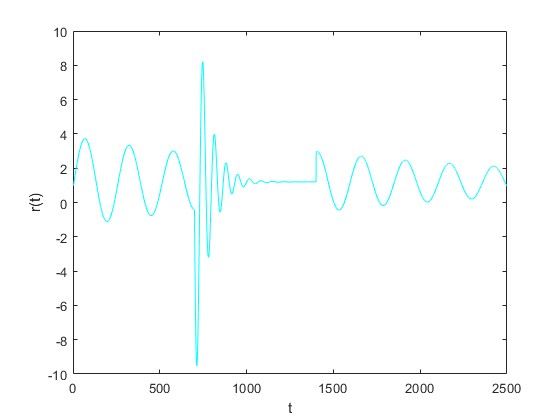}
		\caption{$\alpha=0.86$.}
	\end{subfigure}
	\hfill
	\begin{subfigure}{0.325\textwidth}
		\includegraphics[width=\textwidth]{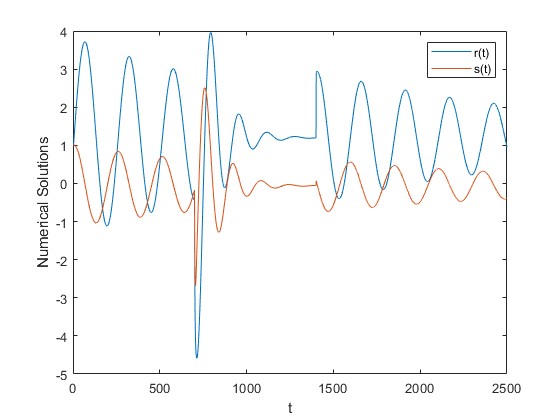}
		\caption{$\alpha=0.89$.}
	\end{subfigure}
	\hfill
	\begin{subfigure}{0.325\textwidth}		\includegraphics[width=\textwidth]{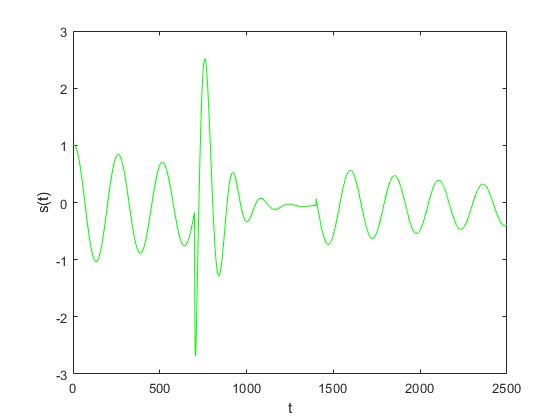}
		\caption{$\alpha=0.89$.}
	\end{subfigure}
	\hfill
	\begin{subfigure}{0.325\textwidth}
		\includegraphics[width=\textwidth]{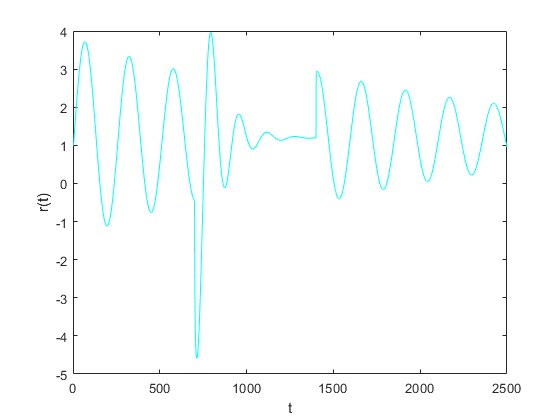}
		\caption{$\alpha=0.89$.}
	\end{subfigure}
	\hfill
	\begin{subfigure}{0.325\textwidth}
		\includegraphics[width=\textwidth]{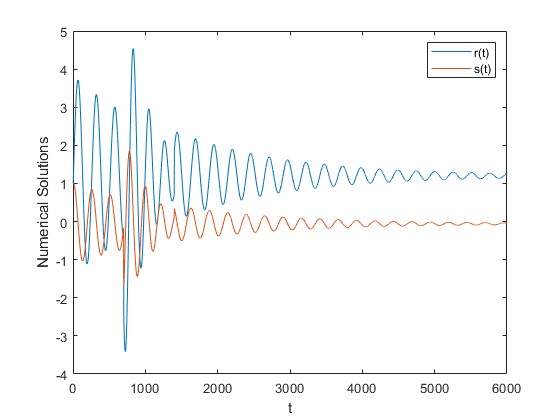}
		\caption{$\alpha=0.95$.}
	\end{subfigure}
	\hfill
	\begin{subfigure}{0.325\textwidth}
		\includegraphics[width=\textwidth]{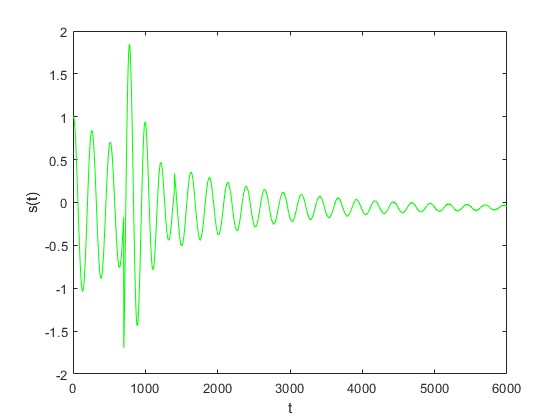}
		\caption{$\alpha=0.95$.}
	\end{subfigure}
	\hfill
	\begin{subfigure}{0.325\textwidth}
		\includegraphics[width=\textwidth]{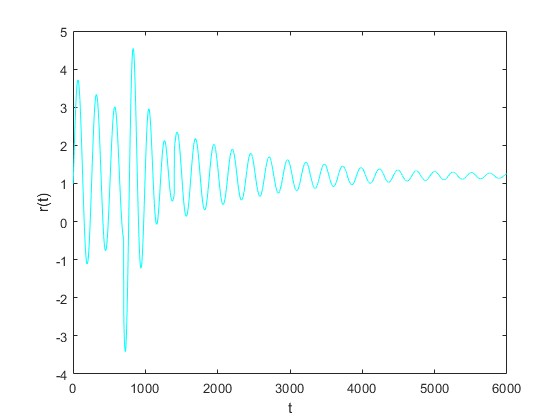}
		\caption{$\alpha=0.95$.}
	\end{subfigure}
	\hfill
	\begin{subfigure}{0.325\textwidth}
		\includegraphics[width=\textwidth]{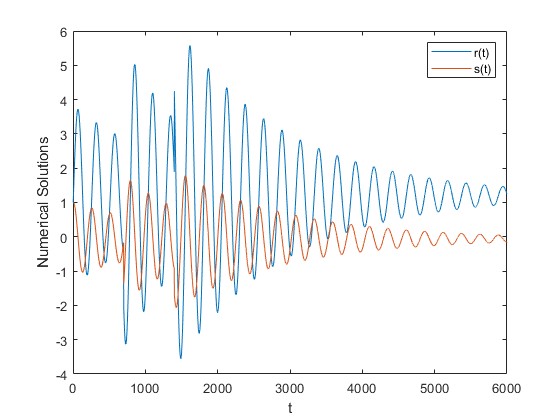}
		\caption{$\alpha=0.99$.}
	\end{subfigure}
	\hfill
	\begin{subfigure}{0.325\textwidth}
		\includegraphics[width=\textwidth]{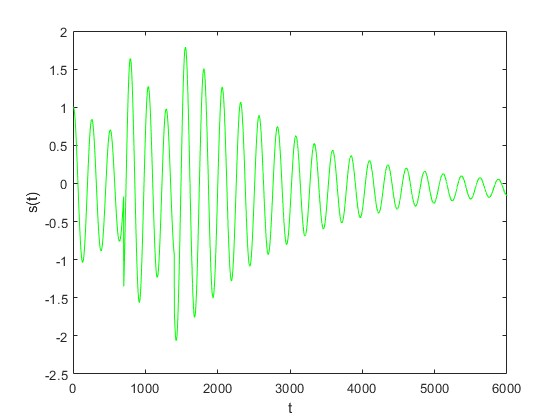}
		\caption{$\alpha=0.99$.}
	\end{subfigure}
	\hfill
	\begin{subfigure}{0.325\textwidth}
		\includegraphics[width=\textwidth]{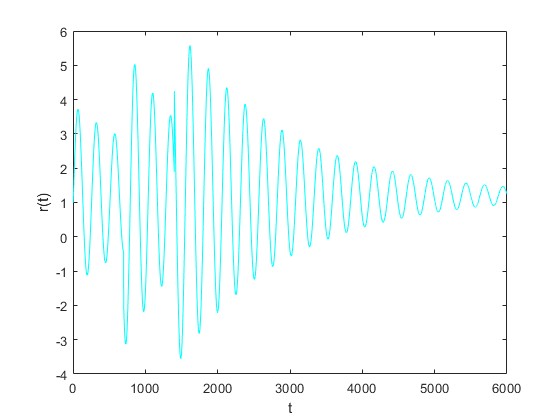}
		\caption{$\alpha=0.99$.}
	\end{subfigure}

	\caption{The Numerical simulations for first model \eqref{ABlinear} with  $t(0)=0$, $s(0)=1$, $p(0)=1$, $\rho_{1}=0.12$, $\rho_{2}=0.05$, $\psi_{1}=0.8$, $\psi_{2}=0.81$, $\gamma_{1}=0.5$, $\gamma_{2}=1.2$, $\sigma_{1}=0.02$, $\sigma_{2}=0.01$, $\omega_{1}=6.1$, $\omega_{2}=-1$.}
	\label{fig:fig5l}
\end{figure}

\begin{figure}[H]
	\centering
	\begin{subfigure}{0.32\textwidth}
		\includegraphics[width=\textwidth]{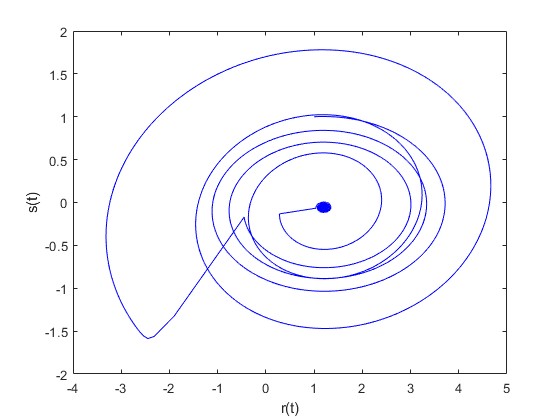}
		\caption{$\rho_{2}=0.01$, $\omega_{2}=-1$.}
	\end{subfigure}
	\hfill
	\begin{subfigure}{0.32\textwidth}
		\includegraphics[width=\textwidth]{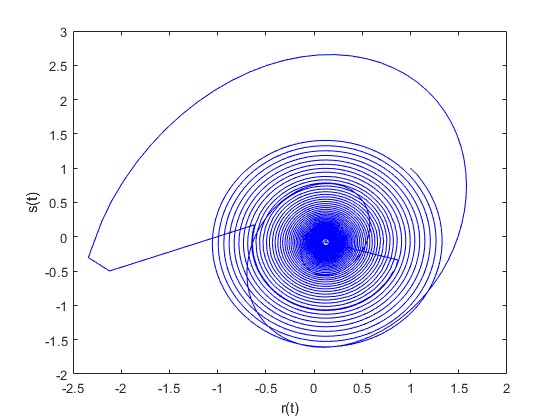}
		\caption{$\rho_{2}=0.01$, $\omega_{2}=-10$.}
	\end{subfigure}
	\hfill
	\begin{subfigure}{0.32\textwidth}
		\includegraphics[width=\textwidth]{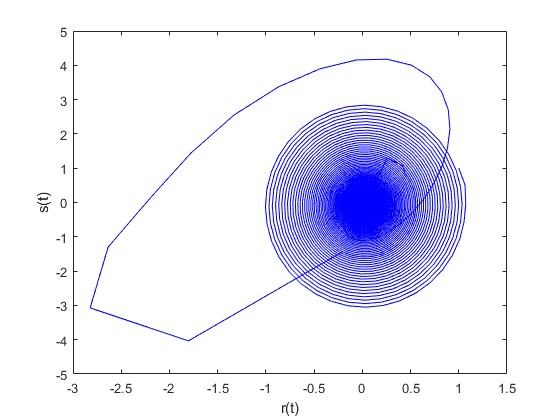}
		\caption{$\rho_{2}=0.01$, $\omega_{2}=-50$.}
	\end{subfigure}
	\hfill
	\begin{subfigure}{0.32\textwidth}
		\includegraphics[width=\textwidth]{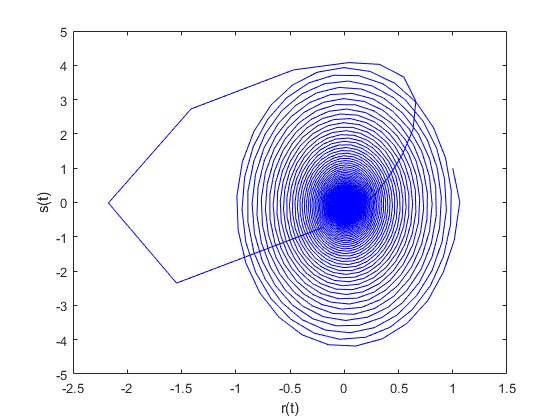}
		\caption{$\rho_{2}=0.01$, $\omega_{2}=-100$.}
	\end{subfigure}
	\hfill
	\begin{subfigure}{0.32\textwidth}
		\includegraphics[width=\textwidth]{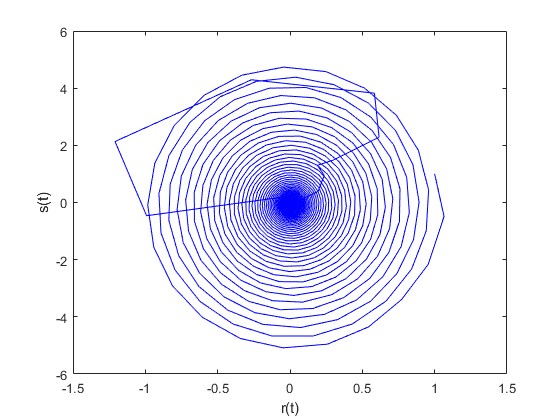}
		\caption{$\rho_{2}=0.01$, $\omega_{2}=-150$.}
	\end{subfigure}
	\hfill
	\begin{subfigure}{0.32\textwidth}
		\includegraphics[width=\textwidth]{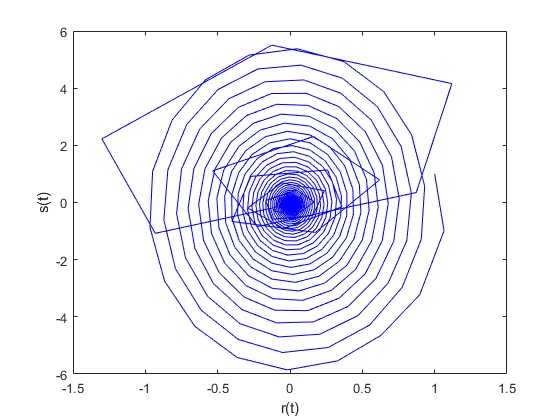}
		\caption{$\rho_{2}=0.01$, $\omega_{2}=-200$.}
	\end{subfigure}
	\hfill
	\begin{subfigure}{0.32\textwidth}
		\includegraphics[width=\textwidth]{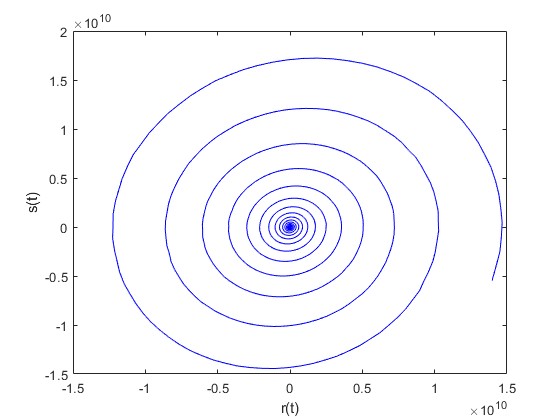}
		\caption{$\rho_{2}=-1$, $\omega_{2}=-10$.}
	\end{subfigure}
	\hfill
	\begin{subfigure}{0.32\textwidth}
		\includegraphics[width=\textwidth]{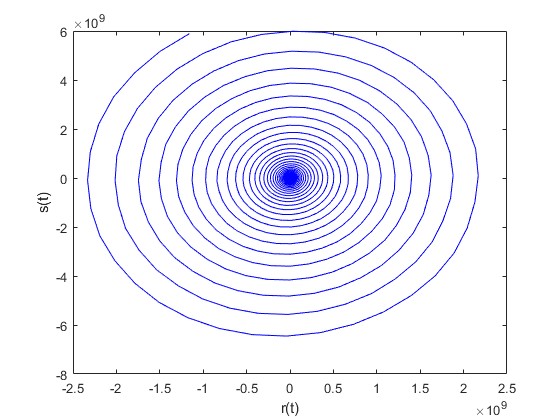}
		\caption{$\rho_{2}=-1$, $\omega_{2}=-50$.}
	\end{subfigure}
	\hfill
	\begin{subfigure}{0.32\textwidth}
		\includegraphics[width=\textwidth]{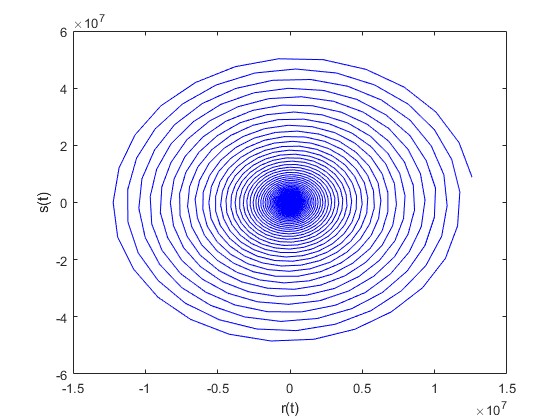}
		\caption{$\rho_{2}=-1$, $\omega_{2}=-100$.}
	\end{subfigure}
	\hfill
	\begin{subfigure}{0.32\textwidth}
		\includegraphics[width=\textwidth]{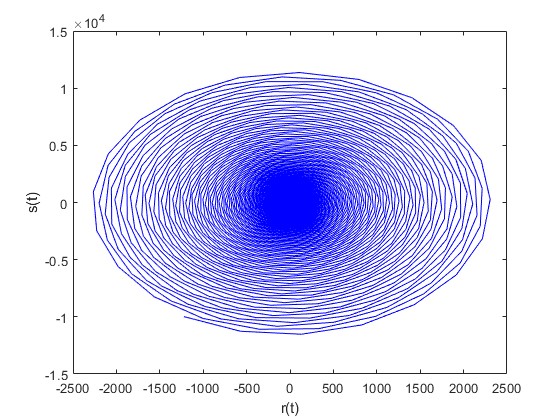}
		\caption{$\rho_{2}=-1$, $\omega_{2}=-100$.}
	\end{subfigure}
	\hfill
	\begin{subfigure}{0.32\textwidth}
		\includegraphics[width=\textwidth]{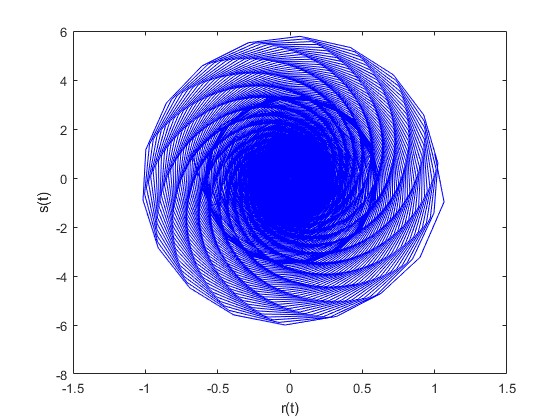}
		\caption{$\rho_{2}=-1$, $\omega_{2}=-200$.}
	\end{subfigure}
	\hfill
	\begin{subfigure}{0.32\textwidth}
		\includegraphics[width=\textwidth]{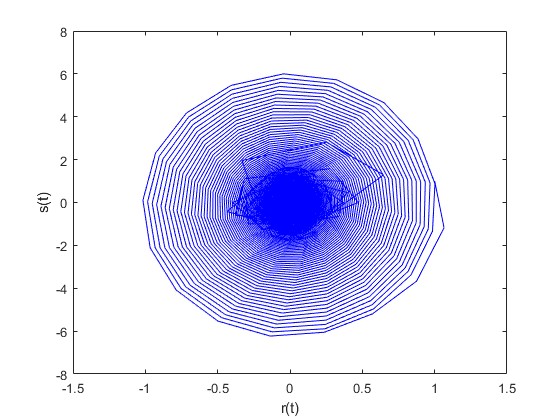}
		\caption{$\rho_{2}=-1$, $\omega_{2}=-225$.}
	\end{subfigure}
	\caption{Chaotic dynamics for the first model \eqref{ABlinear} with $t(0)=0$, $s(0)=1$, $p(0)=1$, $\rho_{1}=0.12$, $\psi_{1}=0.8$, $\psi_{2}=0.81$, $\gamma_{1}=0.5$, $\gamma_{2}=1.2$, $\sigma_{1}=0.02$, $\sigma_{2}=0.01$, $\omega_{1}=6.1$, and $\alpha=0.96$.}
	\label{fig:fig6l}
\end{figure}

\begin{figure}[H]
	\centering
	\begin{subfigure}{0.325\textwidth}
		\includegraphics[width=\textwidth]{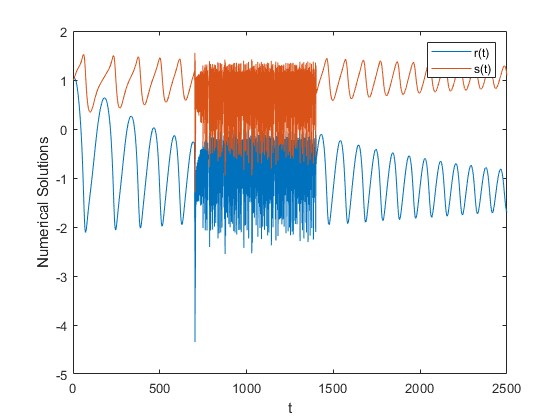}
		\caption{$\alpha=0.81$, $\epsilon=1$}
	\end{subfigure}
	\hfill
	\begin{subfigure}{0.325\textwidth}
		\includegraphics[width=\textwidth]{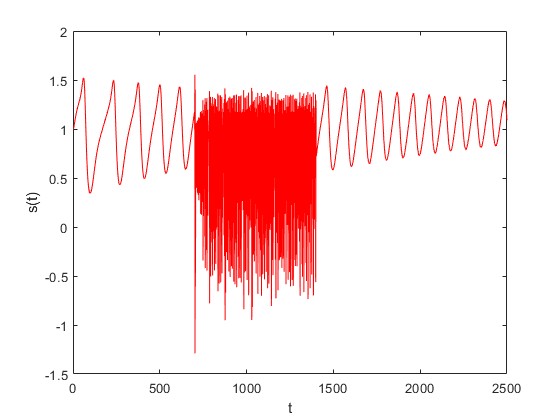}
		\caption{$\alpha=0.81$, $\epsilon=1$.}
	\end{subfigure}
	\hfill
	\begin{subfigure}{0.325\textwidth}
		\includegraphics[width=\textwidth]{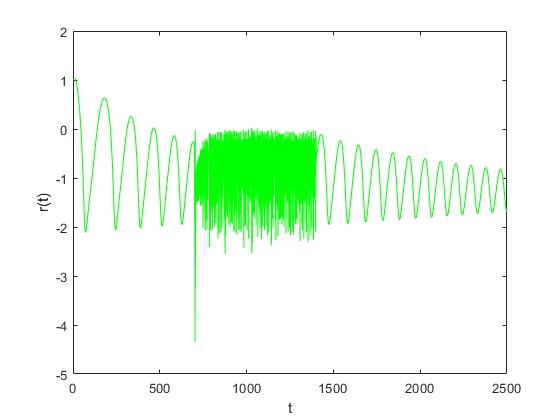}
		\caption{$\alpha=0.81$, $\epsilon=1$.}
	\end{subfigure}
	\hfill
	\begin{subfigure}{0.325\textwidth}
		\includegraphics[width=\textwidth]{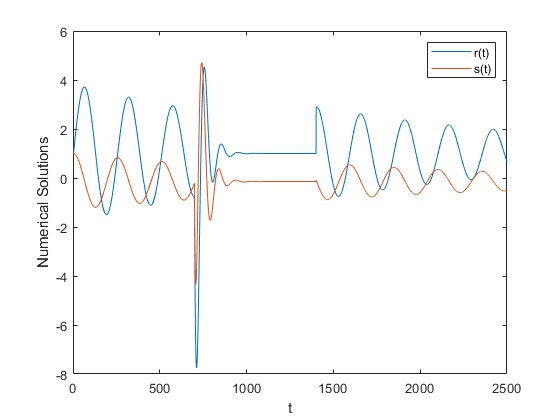}
		\caption{$\alpha=0.87$, $\epsilon=0$.}
	\end{subfigure}
	\hfill
	\begin{subfigure}{0.325\textwidth}		\includegraphics[width=\textwidth]{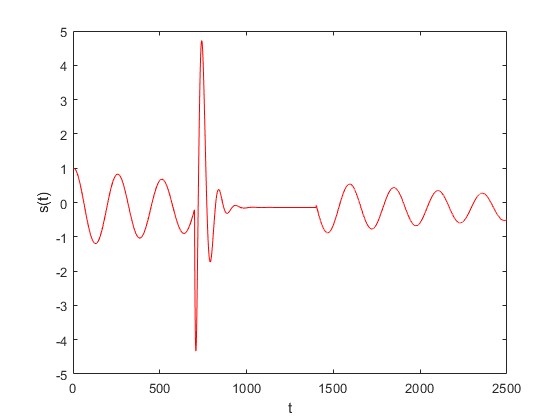}
		\caption{$\alpha=0.87$, $\epsilon=0$.}
	\end{subfigure}
	\hfill
	\begin{subfigure}{0.325\textwidth}
		\includegraphics[width=\textwidth]{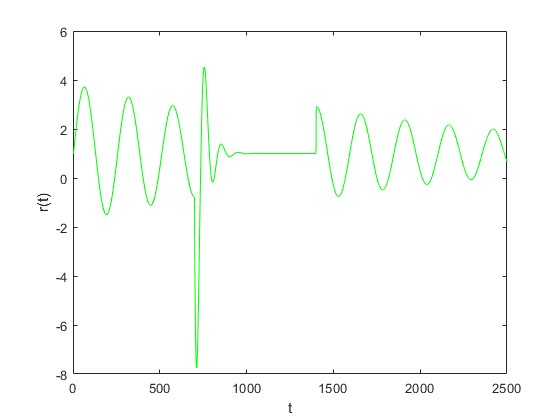}
		\caption{$\alpha=0.87$, $\epsilon=0$.}
	\end{subfigure}
	\hfill
	\begin{subfigure}{0.325\textwidth}
		\includegraphics[width=\textwidth]{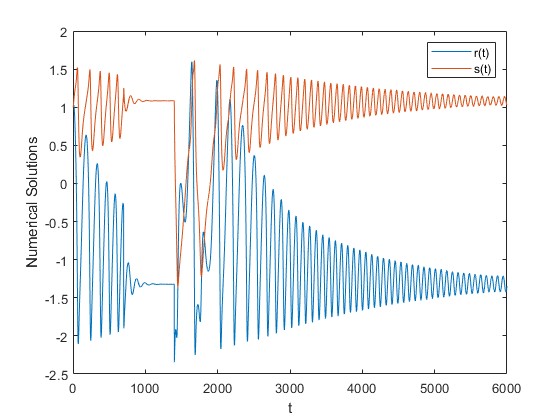}
		\caption{$\alpha=0.94$, $\epsilon=1$.}
	\end{subfigure}
	\hfill
	\begin{subfigure}{0.325\textwidth}
		\includegraphics[width=\textwidth]{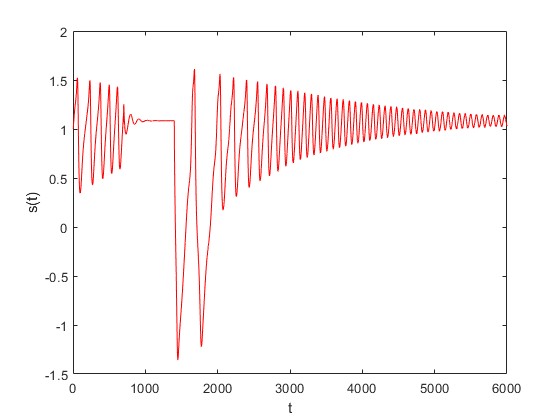}
		\caption{$\alpha=0.94$, $\epsilon=1$.}
	\end{subfigure}
	\hfill
	\begin{subfigure}{0.325\textwidth}
		\includegraphics[width=\textwidth]{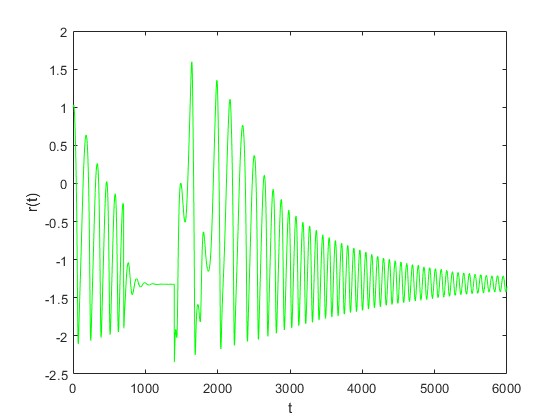}
		\caption{$\alpha=0.94$, $\epsilon=1$.}
	\end{subfigure}
	\hfill
	\begin{subfigure}{0.325\textwidth}
		\includegraphics[width=\textwidth]{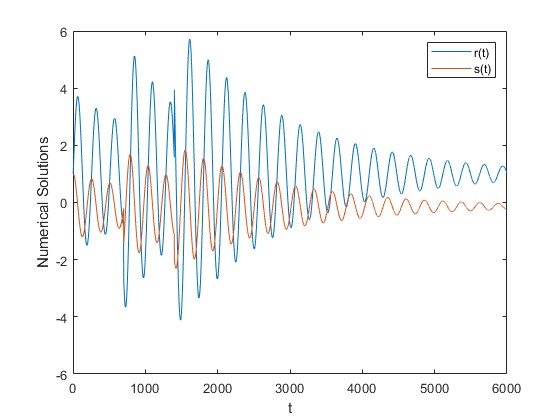}
		\caption{$\alpha=0.99$, $\epsilon=0$.}
	\end{subfigure}
	\hfill
	\begin{subfigure}{0.325\textwidth}
		\includegraphics[width=\textwidth]{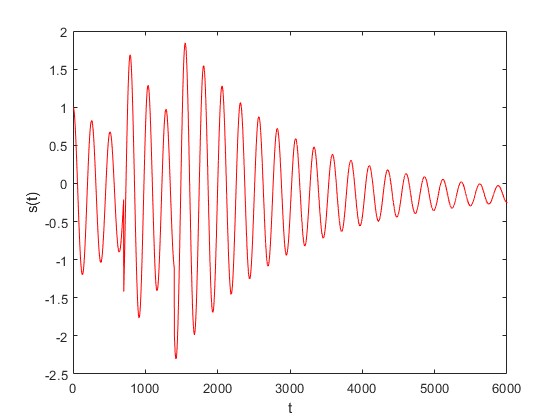}
		\caption{$\alpha=0.99$, $\epsilon=0$.}
	\end{subfigure}
	\hfill
	\begin{subfigure}{0.325\textwidth}
		\includegraphics[width=\textwidth]{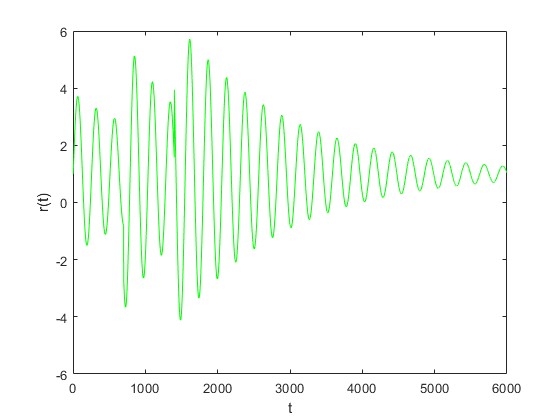}
		\caption{$\alpha=0.99$, $\epsilon=0$.}
	\end{subfigure}

	\caption{The Numerical simulations for second model \eqref{ABnonlinear} with  $t(0)=0$, $r(0)=1$, $s(0)=1$,  $\rho_{1}=0.12$, $\rho_{2}=0.01$,
		$\phi_{1}=1$, $\phi_{2}=1$,
		$\sigma_{1}=0.01$, $\sigma_{2}=0.02$,
		$\omega_{1}=6.1$, $\omega_{2}=-1$.}
	\label{fig:fig7l}
\end{figure}

\begin{figure}[H]
	\centering
	\begin{subfigure}{0.32\textwidth}
		\includegraphics[width=\textwidth]{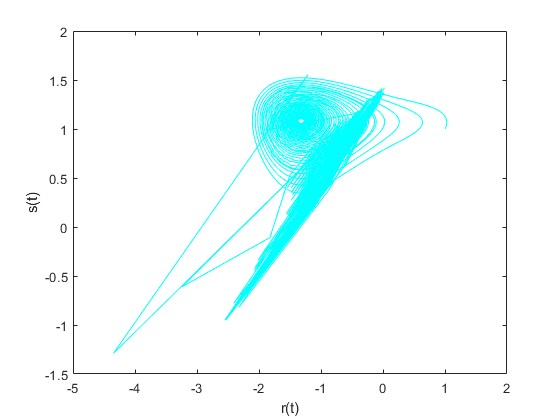}
		\caption{$\alpha=0.81$, $\epsilon=1$.}
	\end{subfigure}
	\hfill
	\begin{subfigure}{0.32\textwidth}
		\includegraphics[width=\textwidth]{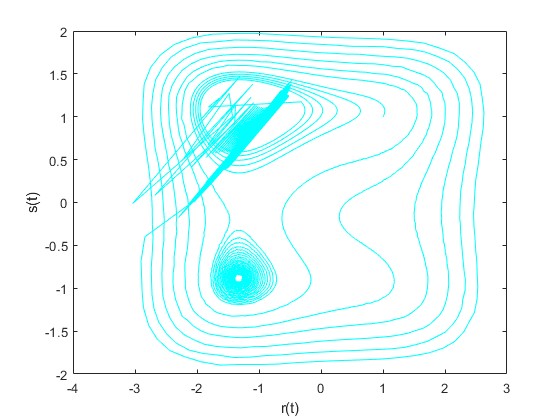}
		\caption{$\alpha=0.87$, $\epsilon=1$.}
	\end{subfigure}
	\hfill
	\begin{subfigure}{0.32\textwidth}
		\includegraphics[width=\textwidth]{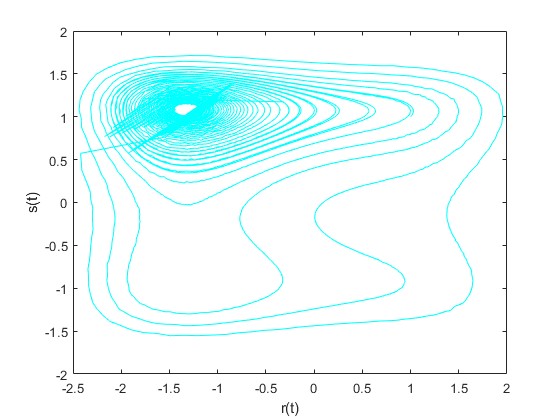}
		\caption{$\alpha=0.92$, $\epsilon=1$.}
	\end{subfigure}
	\hfill
	\begin{subfigure}{0.32\textwidth}
		\includegraphics[width=\textwidth]{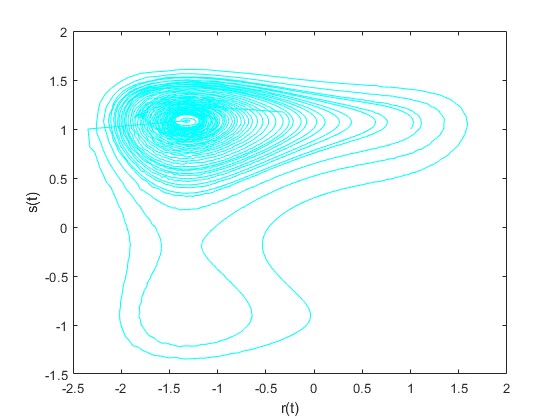}
		\caption{$\alpha=0.95$, $\epsilon=1$.}
	\end{subfigure}
	\hfill
	\begin{subfigure}{0.32\textwidth}
		\includegraphics[width=\textwidth]{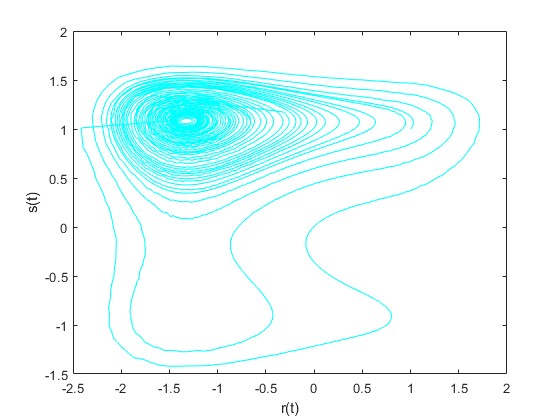}
		\caption{$\alpha=0.99$, $\epsilon=1$.}
	\end{subfigure}
	\hfill
	\begin{subfigure}{0.32\textwidth}
		\includegraphics[width=\textwidth]{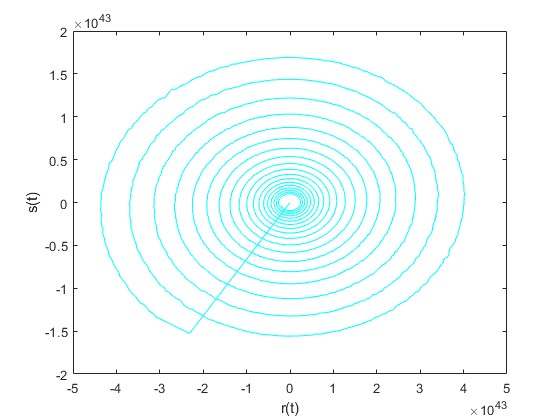}
		\caption{$\alpha=0.81$, $\epsilon=0$.}
	\end{subfigure}
	\hfill
	\begin{subfigure}{0.32\textwidth}
		\includegraphics[width=\textwidth]{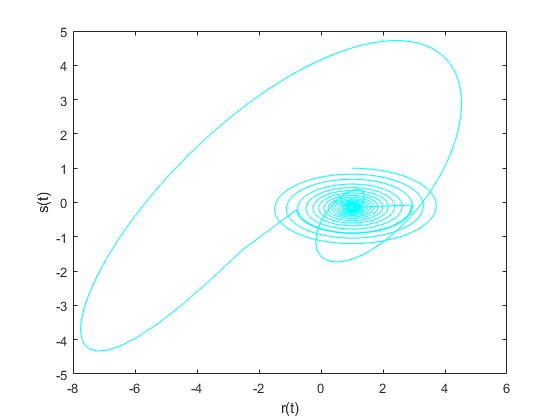}
		\caption{$\alpha=0.87$, $\epsilon=0$.}
	\end{subfigure}
	\hfill
	\begin{subfigure}{0.32\textwidth}
		\includegraphics[width=\textwidth]{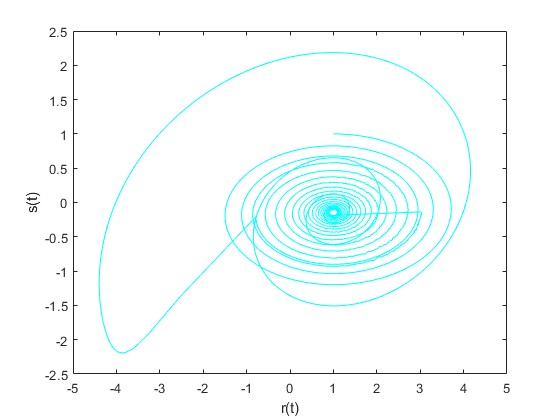}
		\caption{$\alpha=0.92$, $\epsilon=0$.}
	\end{subfigure}
	\hfill
	\begin{subfigure}{0.32\textwidth}
		\includegraphics[width=\textwidth]{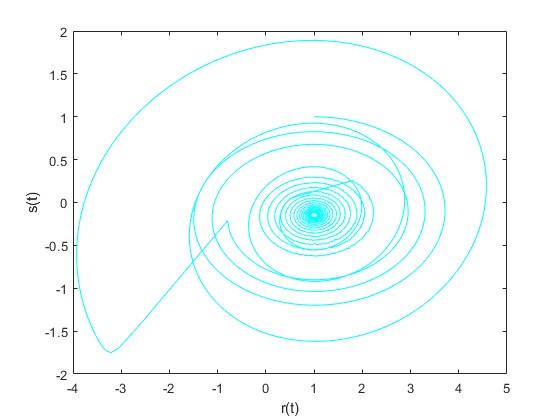}
		\caption{$\alpha=0.95$, $\epsilon=0$.}
	\end{subfigure}
	\hfill
	\begin{subfigure}{0.32\textwidth}
		\includegraphics[width=\textwidth]{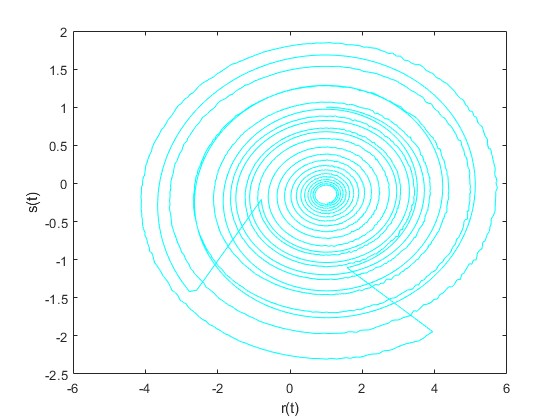}
		\caption{$\alpha=0.99$, $\epsilon=0$.}
	\end{subfigure}
	\hfill
	\begin{subfigure}{0.32\textwidth}
		\includegraphics[width=\textwidth]{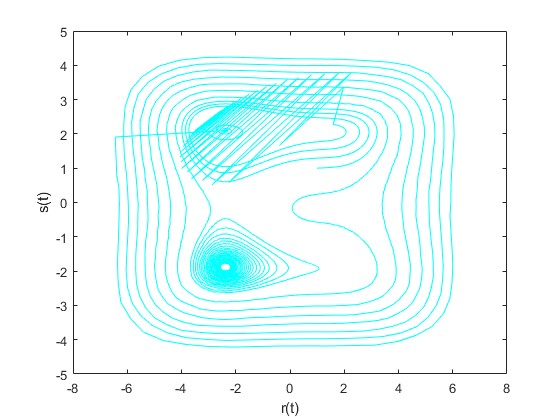}
		\caption{$\alpha=0.95$, $\epsilon=0.25$.}
	\end{subfigure}
	\hfill
	\begin{subfigure}{0.32\textwidth}
		\includegraphics[width=\textwidth]{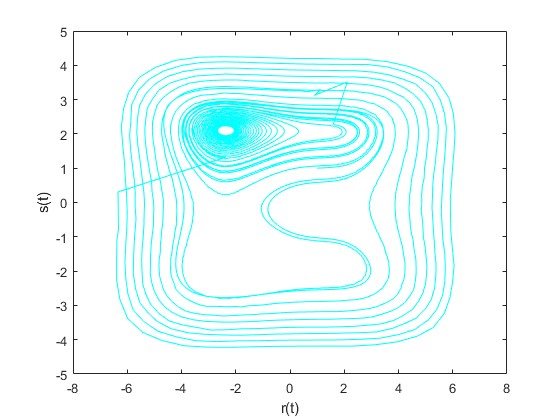}
		\caption{$\alpha=0.99$, $\epsilon=0.25$.}
	\end{subfigure}
	\caption{Chaotic dynamics for the second model \eqref{ABnonlinear} with  $t(0)=0$, $r(0)=1$, $s(0)=1$,  $\rho_{1}=0.12$, $\rho_{2}=0.01$,
		$\phi_{1}=1$, $\phi_{2}=1$,
		$\sigma_{1}=0.01$, $\sigma_{2}=0.02$,
		$\omega_{1}=6.1$, $\omega_{2}=-1$.}
	\label{fig:fig8l}
\end{figure}

	It presents the numerical simulations and chaotic phase portrait via first \eqref{ABlinear} and second \eqref{ABnonlinear}models for love dynamical system.
	The figures \ref{fig:fig5l}(a)-\ref{fig:fig5l}(l), and \ref{fig:fig6l}(a)-\ref{fig:fig6l}(l) demonstrate the numerical simulations, while the figures \ref{fig:fig7l}(a)-\ref{fig:fig7l}(l), and \ref{fig:fig8l}(a)-\ref{fig:fig8l}(l) demonstrate the chaotic behaviors of phase portrait for the system under consideration.

\subsection{Numerical algorithm in the aspect of Caputo-Fabrizio derivative:}
We now investigate the third case \eqref{linear_caputo_fading}. The numerical scheme is obtained using Newton interpolation \eqref{Newton_method} and is given as
\begin{align}
	\begin{cases}
		&U^{k_{1}} = U(0)+\frac{1}{12}\sum_{j_1 = 2}^{k_{1}}\bigg[23e(t_{j_1}, U(t_{j_1}))-16e(t_{j_1 - 1}, U(t_{j_1 - 1}))+5e(t_{j_1 - 2}, U(t_{j_1 - 2})) \bigg]* \Delta t, 0 \le t\le a_{1}  \\
		&U^{k_{2}}= U(a_{1})+\frac{1-\gamma}{M(\gamma)} \bigg[ \psi(t_{k}, U(t_{k})) - e(t_{k - 1}, U(t_{k - 1})) \bigg] + \\ &~\frac{1}{12}\sum_{j_2 = k_{1}+1}^{k_{2}} \bigg[23e(t_{j_2}, U(t_{j_2})) -16e(t_{j_2 - 1}, U(t_{j_2 - 1}))+5e(t_{j_2 - 2}, U(t_{j_2 - 2})) \bigg]*\Delta t, 	a_{1} \le t\le a_2\\
		&U^{k_{3}} = U(a_{2})+ \frac{1}{12} \sum_{j_3 = k_{2}+3}^{k_{3}}\bigg[23e(t_{j_3}, U(t_{j_3}))- 16e(t_{j_3 - 1}, U(t_{j_3 - 1})) + 5e(t_{j_3 - 2}, U(t_{j_3 - 2})) \bigg]*\Delta t+ \\
		&~\sigma_{i} \sum_{j_3 = k_{2}+3}^{k_{3}}U(\mathbb{B}_i(t_{j_3})-\mathbb{B}_i(t_{j_3-1})), a_{2} \le t\le a 
	\end{cases}
\end{align}
\subsection*{ Computer simulation of first and second modified love dynamical systems in the aspect of Caputo-Fabrizio derivative:}  
\begin{figure}[H]
	\centering
	\begin{subfigure}{0.325\textwidth}
		\includegraphics[width=\textwidth]{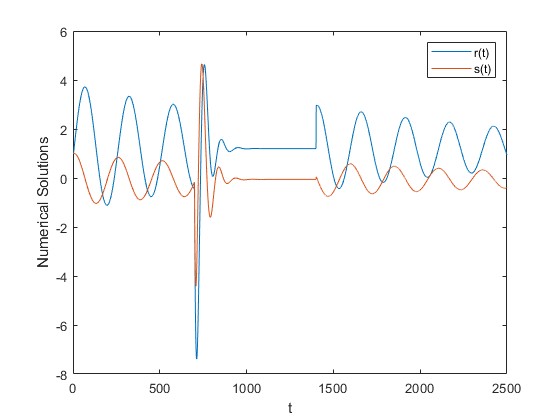}
		\caption{$\alpha=0.87$.}
	\end{subfigure}
	\hfill
	\begin{subfigure}{0.325\textwidth}
		\includegraphics[width=\textwidth]{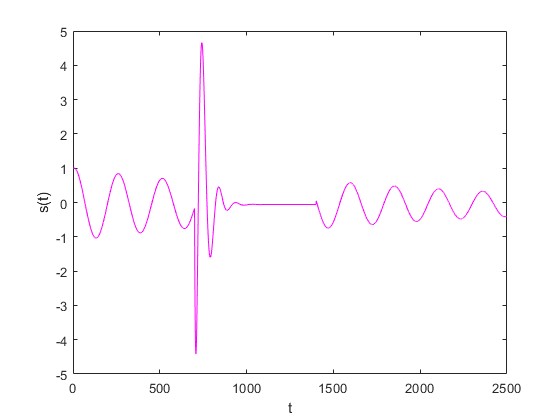}
		\caption{$\alpha=0.87$.}
	\end{subfigure}
	\hfill
	\begin{subfigure}{0.325\textwidth}
		\includegraphics[width=\textwidth]{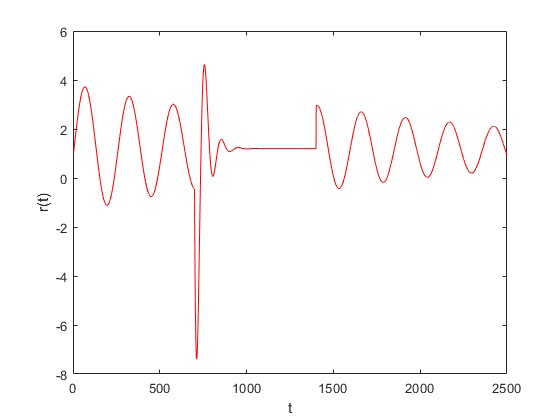}
		\caption{$\alpha=0.87$.}
	\end{subfigure}
	\hfill
	\begin{subfigure}{0.325\textwidth}
		\includegraphics[width=\textwidth]{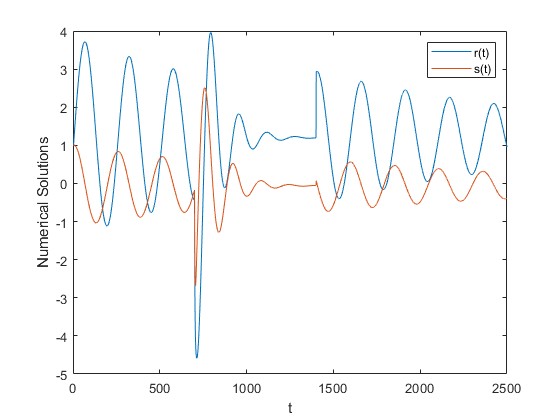}
		\caption{$\alpha=0.9$.}
	\end{subfigure}
	\hfill
	\begin{subfigure}{0.325\textwidth}		\includegraphics[width=\textwidth]{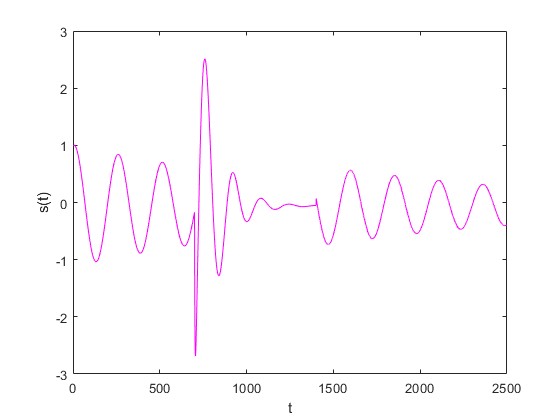}
		\caption{$\alpha=0.9$.}
	\end{subfigure}
	\hfill
	\begin{subfigure}{0.325\textwidth}
		\includegraphics[width=\textwidth]{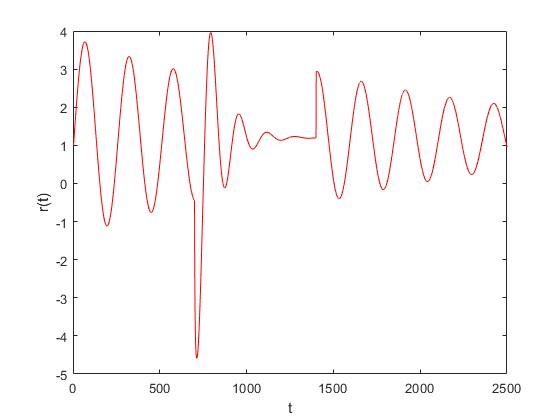}
		\caption{$\alpha=0.9$.}
	\end{subfigure}
	\hfill
	\begin{subfigure}{0.325\textwidth}
		\includegraphics[width=\textwidth]{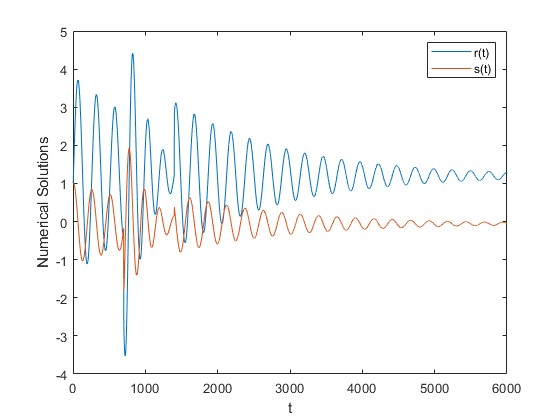}
		\caption{$\alpha=0.94$.}
	\end{subfigure}
	\hfill
	\begin{subfigure}{0.325\textwidth}
		\includegraphics[width=\textwidth]{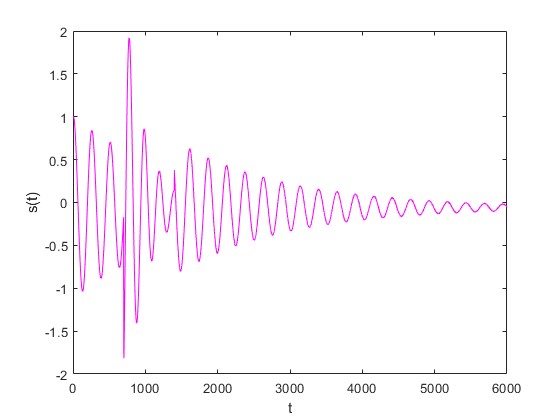}
		\caption{$\alpha=0.94$.}
	\end{subfigure}
	\hfill
	\begin{subfigure}{0.325\textwidth}
		\includegraphics[width=\textwidth]{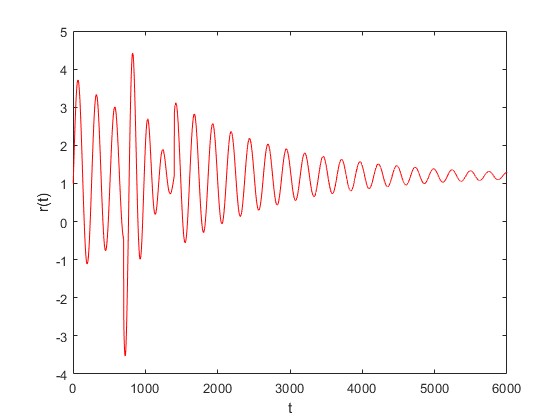}
		\caption{$\alpha=0.94$.}
	\end{subfigure}
	\hfill
	\begin{subfigure}{0.325\textwidth}
		\includegraphics[width=\textwidth]{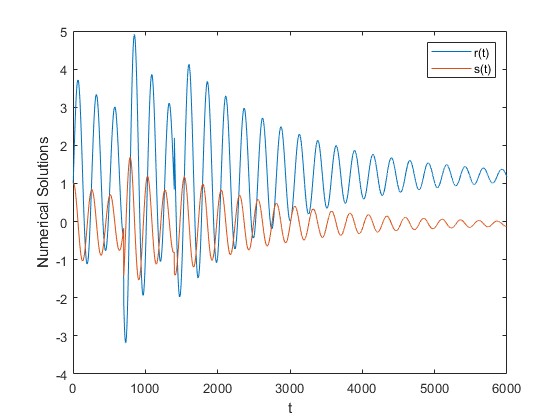}
		\caption{$\alpha=0.98$.}
	\end{subfigure}
	\hfill
	\begin{subfigure}{0.325\textwidth}
		\includegraphics[width=\textwidth]{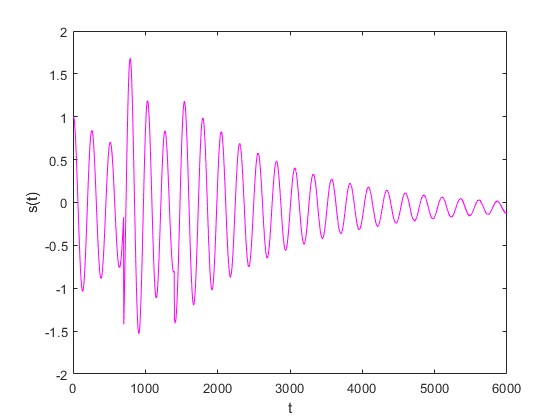}
		\caption{$\alpha=0.98$.}
	\end{subfigure}
	\hfill
	\begin{subfigure}{0.325\textwidth}
		\includegraphics[width=\textwidth]{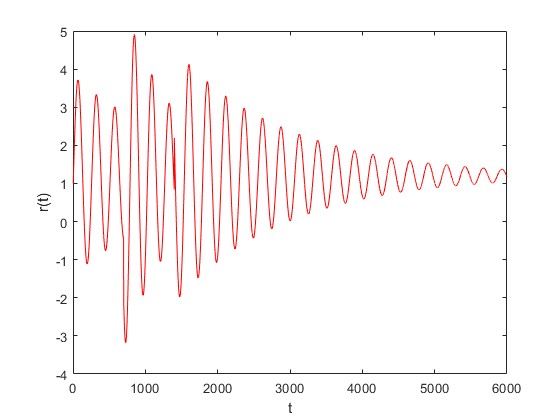}
		\caption{$\alpha=0.98$.}
	\end{subfigure}

	\caption{The Numerical simulations for first model \eqref{CFlinear} with  $t(0)=0$, $s(0)=1$, $p(0)=1$, $\rho_{1}=0.12$, $\rho_{2}=0.05$, $\psi_{1}=0.8$, $\psi_{2}=0.81$, $\gamma_{1}=0.5$, $\gamma_{2}=1.2$, $\sigma_{1}=0.02$, $\sigma_{2}=0.01$, $\omega_{1}=6.1$, $\omega_{2}=-1$. }
	\label{fig:fig9l}
\end{figure}

\begin{figure}[H]
	\centering
	\begin{subfigure}{0.32\textwidth}
		\includegraphics[width=\textwidth]{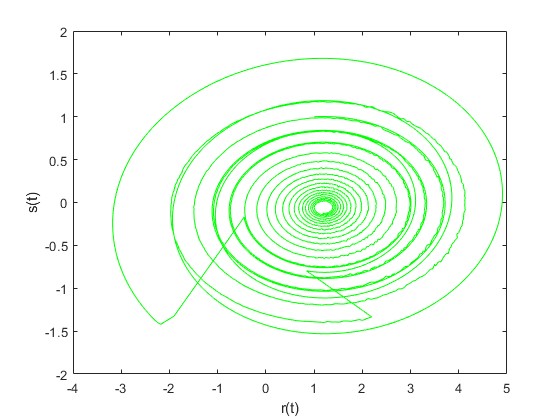}
		\caption{$\rho_{2}=0.01$, $\omega_{2}=-1$}
	\end{subfigure}
	\hfill
	\begin{subfigure}{0.32\textwidth}
		\includegraphics[width=\textwidth]{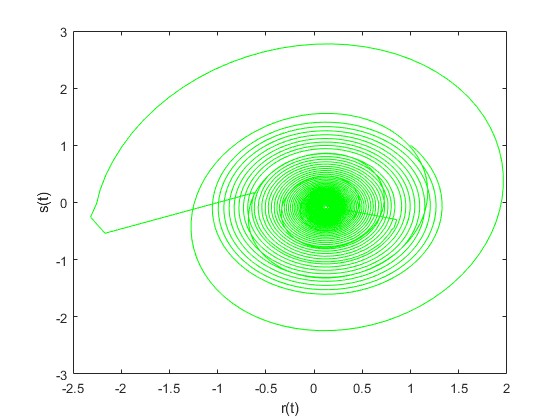}
		\caption{$\rho_{2}=0.01$, $\omega_{2}=-10$}
	\end{subfigure}
	\hfill
	\begin{subfigure}{0.32\textwidth}
		\includegraphics[width=\textwidth]{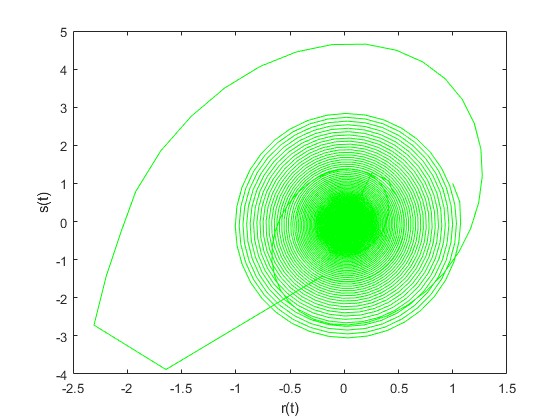}
		\caption{$\rho_{2}=0.01$, $\omega_{2}=-50$}
	\end{subfigure}
	\hfill
	\begin{subfigure}{0.32\textwidth}
		\includegraphics[width=\textwidth]{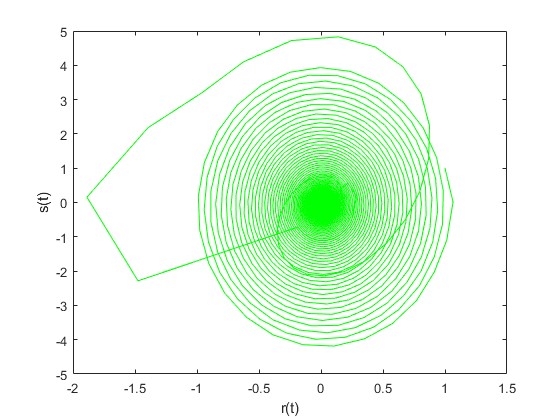}
		\caption{$\rho_{2}=0.01$, $\omega_{2}=-100$}
	\end{subfigure}
	\hfill
	\begin{subfigure}{0.32\textwidth}
		\includegraphics[width=\textwidth]{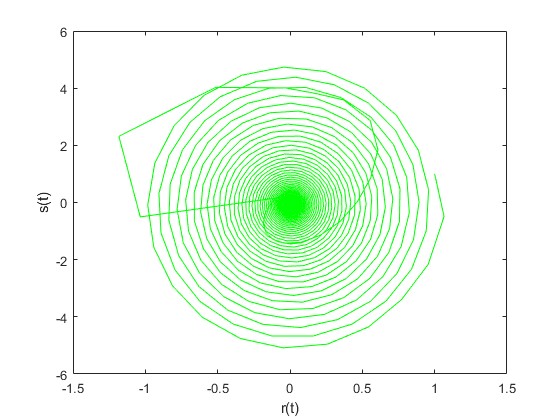}
		\caption{$\rho_{2}=0.01$, $\omega_{2}=-150$}
	\end{subfigure}
	\hfill
	\begin{subfigure}{0.32\textwidth}
		\includegraphics[width=\textwidth]{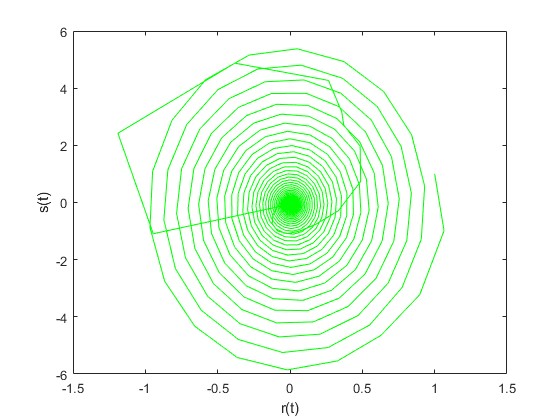}
		\caption{$\rho_{2}=0.01$, $\omega_{2}=-200$}
	\end{subfigure}
	\hfill
	\begin{subfigure}{0.32\textwidth}
		\includegraphics[width=\textwidth]{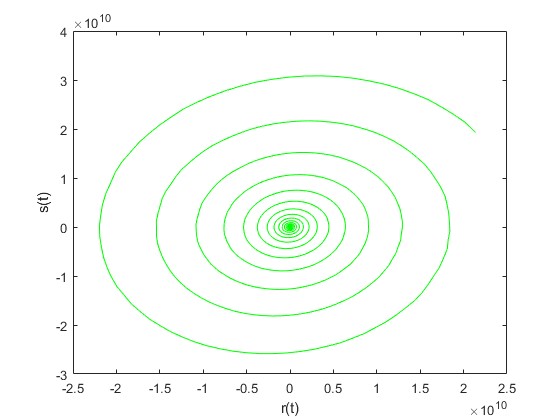}
		\caption{$\rho_{2}=-1$, $\omega_{2}=-10$}
	\end{subfigure}
	\hfill
	\begin{subfigure}{0.32\textwidth}
		\includegraphics[width=\textwidth]{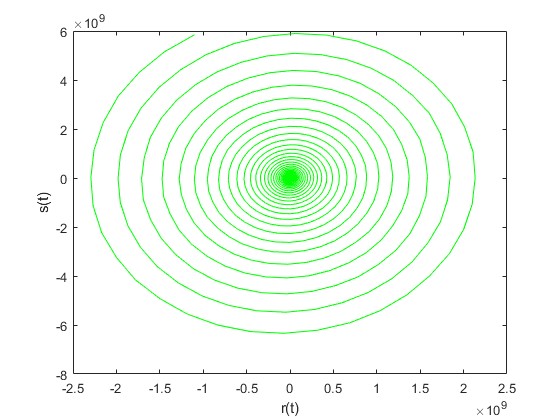}
		\caption{$\rho_{2}=-1$, $\omega_{2}=-50$}
	\end{subfigure}
	\hfill
	\begin{subfigure}{0.32\textwidth}
		\includegraphics[width=\textwidth]{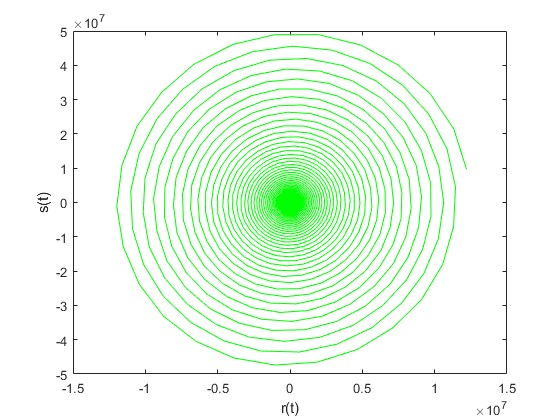}
		\caption{$\rho_{2}=-1$, $\omega_{2}=-100$}
	\end{subfigure}
	\hfill
	\begin{subfigure}{0.32\textwidth}
		\includegraphics[width=\textwidth]{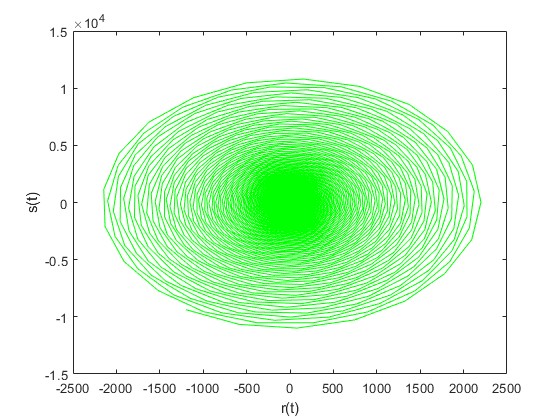}
		\caption{$\rho_{2}=-1$, $\omega_{2}=-150$}
	\end{subfigure}
	\hfill
	\begin{subfigure}{0.32\textwidth}
		\includegraphics[width=\textwidth]{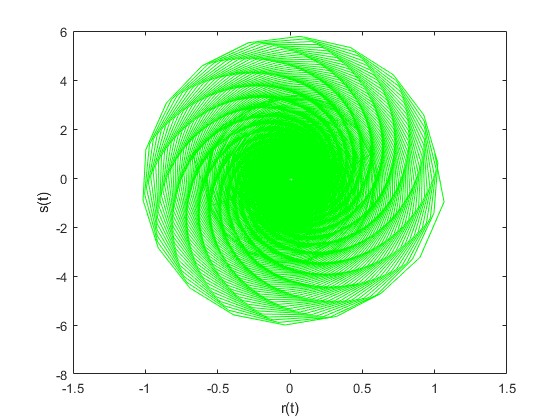}
		\caption{$\rho_{2}=-1$, $\omega_{2}=-200$}
	\end{subfigure}
	\hfill
	\begin{subfigure}{0.32\textwidth}
		\includegraphics[width=\textwidth]{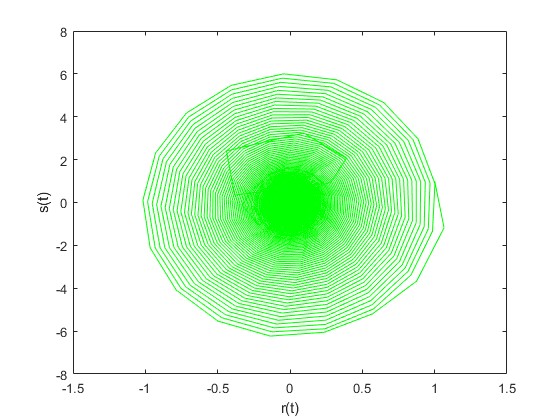}
		\caption{$\rho_{2}=-1$, $\omega_{2}=-225$}
	\end{subfigure}
	\caption{Chaotic dynamics for the first model \eqref{CFlinear} with  $t(0)=0$, $s(0)=1$, $p(0)=1$, $\rho_{1}=0.12$, $\psi_{1}=0.8$, $\psi_{2}=0.81$, $\gamma_{1}=0.5$, $\gamma_{2}=1.2$, $\sigma_{1}=0.02$, $\sigma_{2}=0.01$, $\omega_{1}=6.1$, and $\alpha=0.97$.}
	\label{fig:fig10l}
\end{figure}

\begin{figure}[H]
	\centering
	\begin{subfigure}{0.325\textwidth}
		\includegraphics[width=\textwidth]{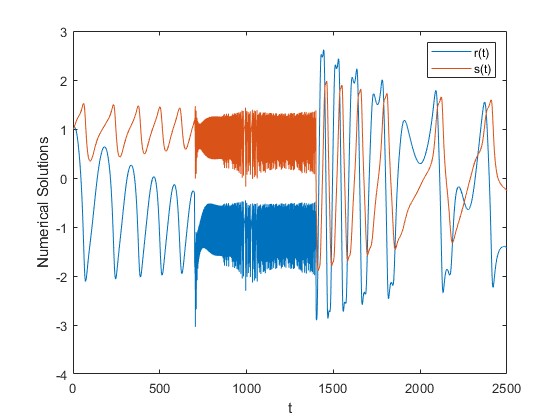}
		\caption{$\alpha=0.81$, $\epsilon=1$}
	\end{subfigure}
	\hfill
	\begin{subfigure}{0.325\textwidth}
		\includegraphics[width=\textwidth]{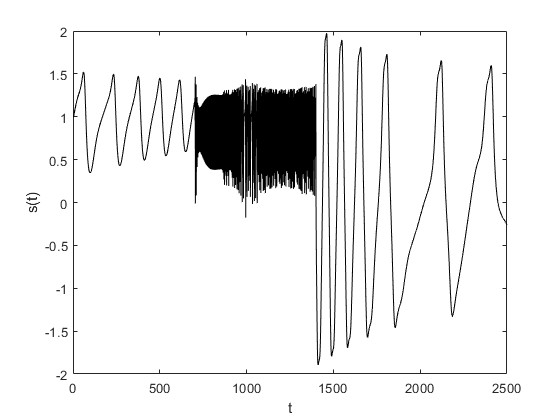}
		\caption{$\alpha=0.81$, $\epsilon=1$.}
	\end{subfigure}
	\hfill
	\begin{subfigure}{0.325\textwidth}
		\includegraphics[width=\textwidth]{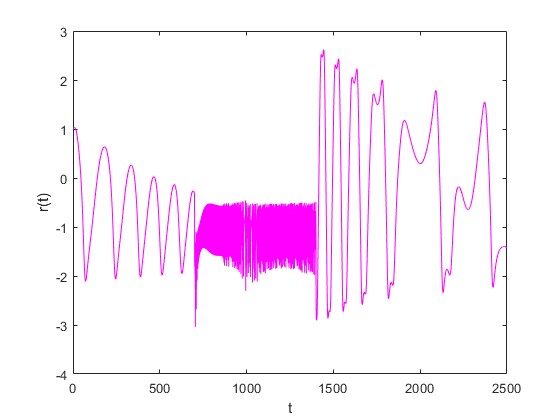}
		\caption{$\alpha=0.81$, $\epsilon=1$.}
	\end{subfigure}
	\hfill
	\begin{subfigure}{0.325\textwidth}
		\includegraphics[width=\textwidth]{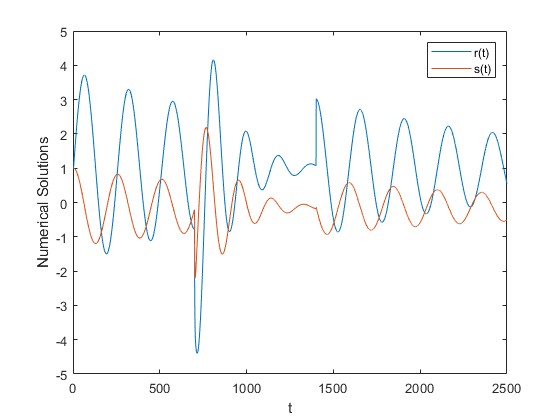}
		\caption{$\alpha=0.87$, $\epsilon=0$.}
	\end{subfigure}
	\hfill
	\begin{subfigure}{0.325\textwidth}		\includegraphics[width=\textwidth]{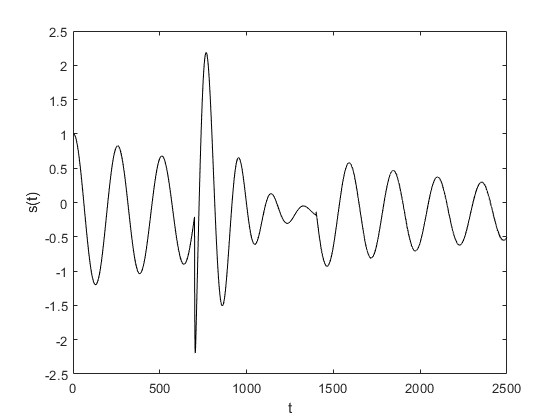}
		\caption{$\alpha=0.87$, $\epsilon=0$.}
	\end{subfigure}
	\hfill
	\begin{subfigure}{0.325\textwidth}
		\includegraphics[width=\textwidth]{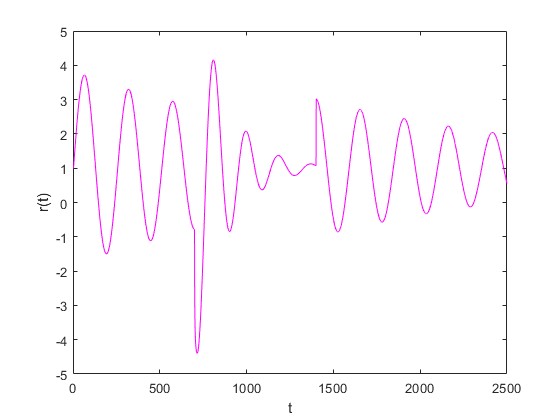}
		\caption{$\alpha=0.87$, $\epsilon=0$.}
	\end{subfigure}
	\hfill
	\begin{subfigure}{0.325\textwidth}
		\includegraphics[width=\textwidth]{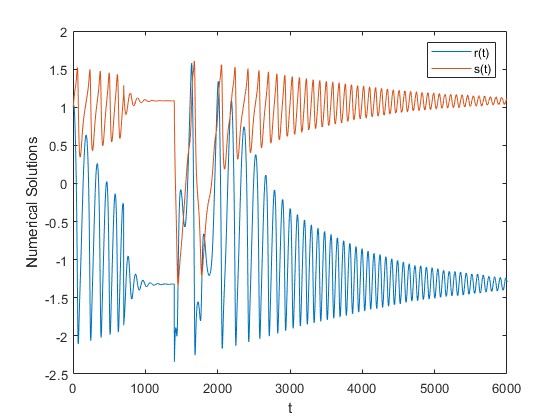}
		\caption{$\alpha=0.94$, $\epsilon=1$.}
	\end{subfigure}
	\hfill
	\begin{subfigure}{0.325\textwidth}
		\includegraphics[width=\textwidth]{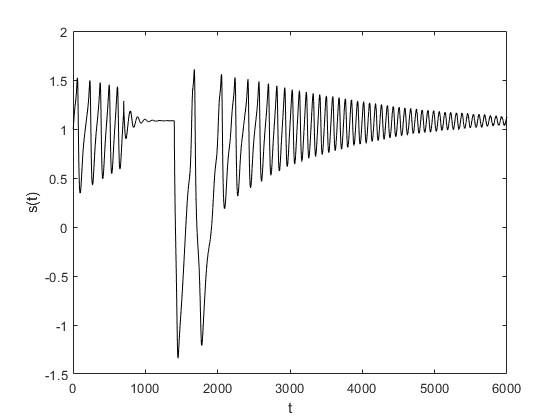}
		\caption{$\alpha=0.94$, $\epsilon=1$.}
	\end{subfigure}
	\hfill
	\begin{subfigure}{0.325\textwidth}
		\includegraphics[width=\textwidth]{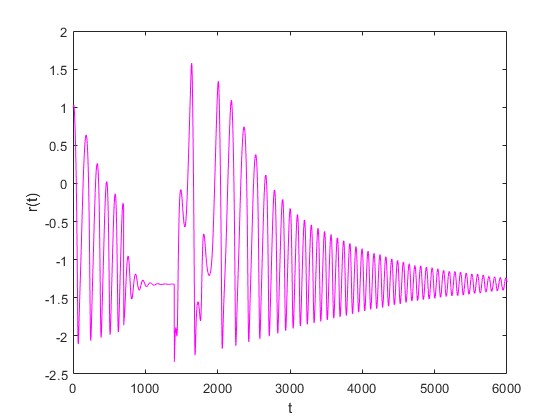}
		\caption{$\alpha=0.94$, $\epsilon=1$.}
	\end{subfigure}
	\hfill
	\begin{subfigure}{0.325\textwidth}
		\includegraphics[width=\textwidth]{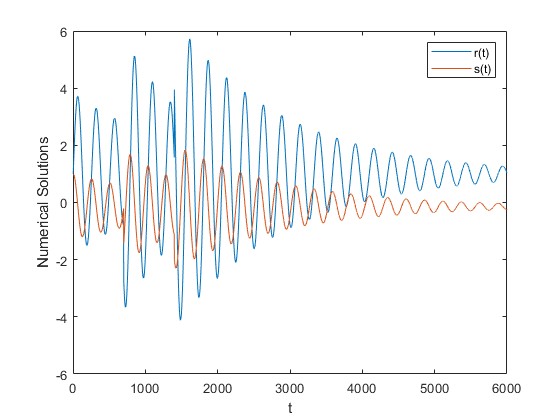}
		\caption{$\alpha=0.99$, $\epsilon=0$.}
	\end{subfigure}
	\hfill
	\begin{subfigure}{0.325\textwidth}
		\includegraphics[width=\textwidth]{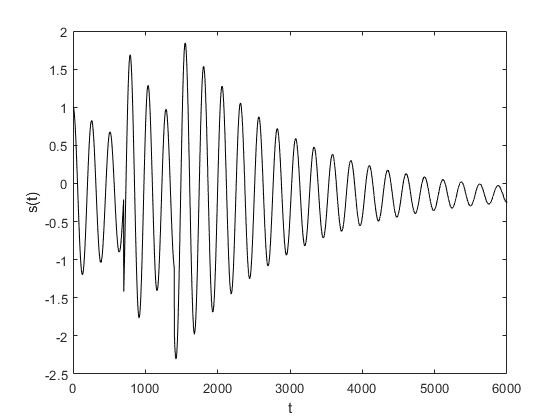}
		\caption{$\alpha=0.99$, $\epsilon=0$.}
	\end{subfigure}
	\hfill
	\begin{subfigure}{0.325\textwidth}
		\includegraphics[width=\textwidth]{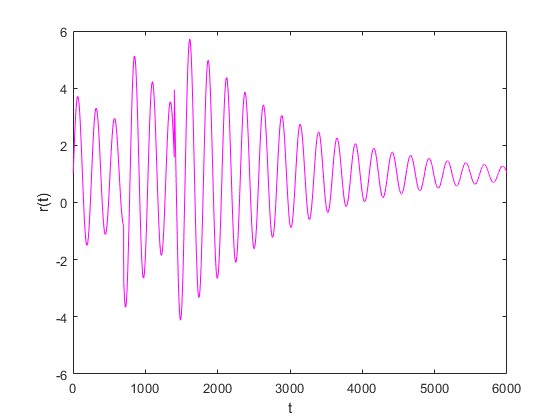}
		\caption{$\alpha=0.99$, $\epsilon=0$.}
	\end{subfigure}

	\caption{The Numerical simulations for second model \eqref{CFnonlinear} with  $t(0)=0$, $r(0)=1$, $s(0)=1$,  $\rho_{1}=0.12$, $\rho_{2}=0.01$,
		$\phi_{1}=1$, $\phi_{2}=1$,
		$\sigma_{1}=0.01$, $\sigma_{2}=0.02$,
		$\omega_{1}=6.1$, $\omega_{2}=-1$.}
	\label{fig:fig11l}
\end{figure}

\begin{figure}[H]
	\centering
	\begin{subfigure}{0.32\textwidth}
		\includegraphics[width=\textwidth]{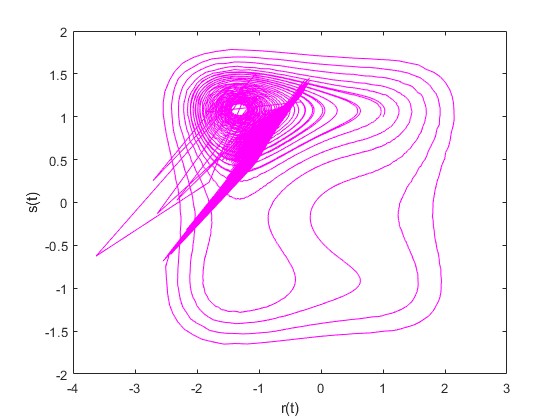}
		\caption{$\alpha=0.84$, $\epsilon=1$.}
	\end{subfigure}
	\hfill
	\begin{subfigure}{0.32\textwidth}
		\includegraphics[width=\textwidth]{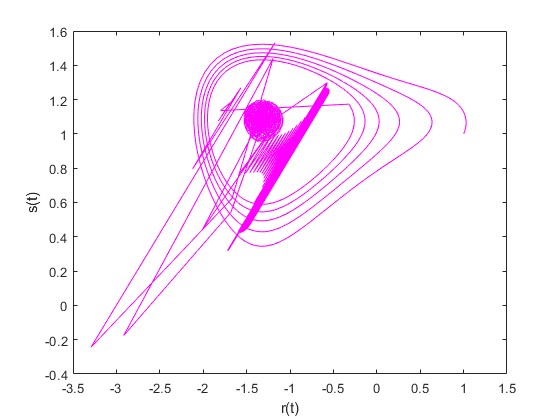}
		\caption{$\alpha=0.88$, $\epsilon=1$.}
	\end{subfigure}
	\hfill
	\begin{subfigure}{0.32\textwidth}
		\includegraphics[width=\textwidth]{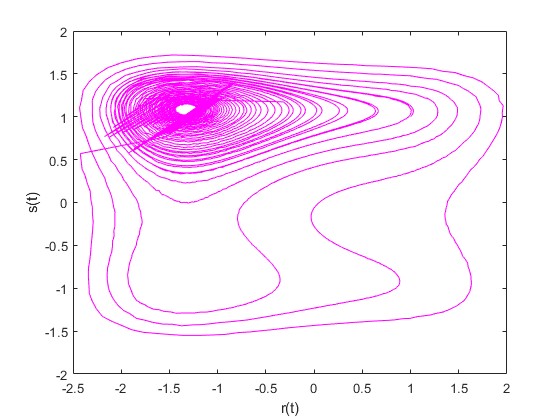}
		\caption{$\alpha=0.92$, $\epsilon=1$.}
	\end{subfigure}
	\hfill
	\begin{subfigure}{0.32\textwidth}
		\includegraphics[width=\textwidth]{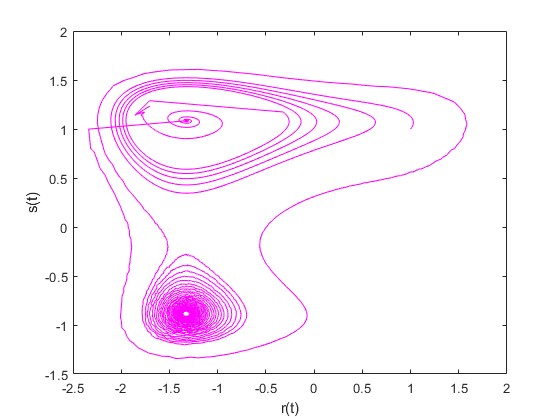}
		\caption{$\alpha=0.96$, $\epsilon=1$.}
	\end{subfigure}
	\hfill
	\begin{subfigure}{0.32\textwidth}
		\includegraphics[width=\textwidth]{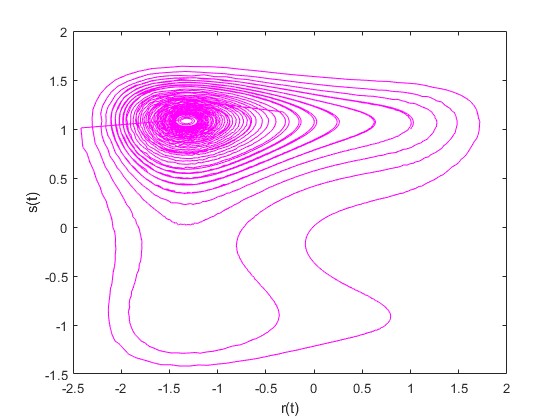}
		\caption{$\alpha=0.99$, $\epsilon=1$.}
	\end{subfigure}
	\hfill
	\begin{subfigure}{0.32\textwidth}
		\includegraphics[width=\textwidth]{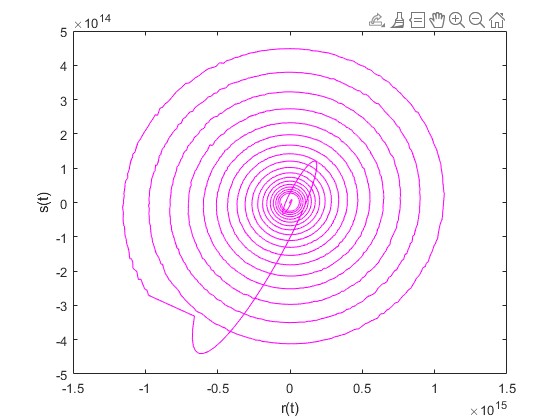}
		\caption{$\alpha=0.84$, $\epsilon=0$.}
	\end{subfigure}
	\hfill
	\begin{subfigure}{0.32\textwidth}
		\includegraphics[width=\textwidth]{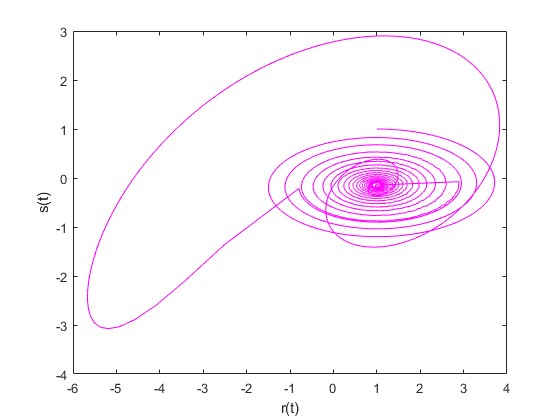}
		\caption{$\alpha=0.89$, $\epsilon=0$.}
	\end{subfigure}
	\hfill
	\begin{subfigure}{0.32\textwidth}
		\includegraphics[width=\textwidth]{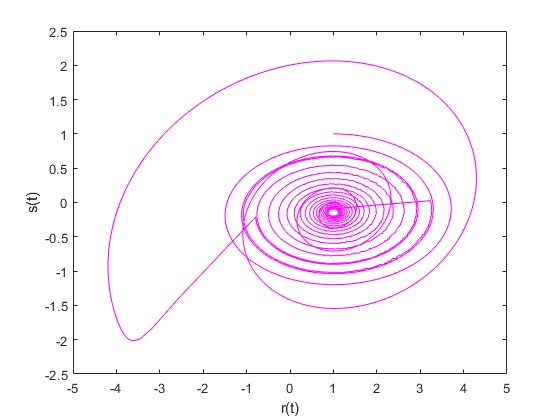}
		\caption{$\alpha=0.93$, $\epsilon=0$.}
	\end{subfigure}
	\hfill
	\begin{subfigure}{0.32\textwidth}
		\includegraphics[width=\textwidth]{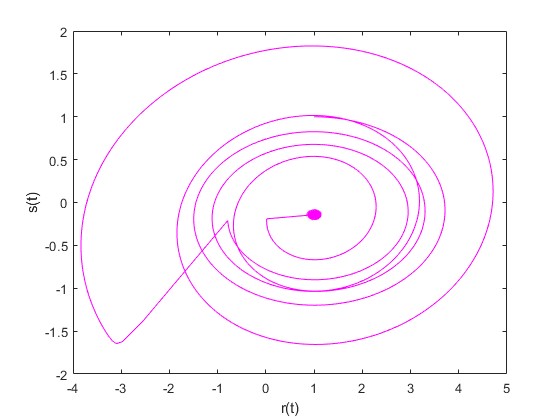}
		\caption{$\alpha=0.96$, $\epsilon=0$.}
	\end{subfigure}
	\hfill
	\begin{subfigure}{0.32\textwidth}
		\includegraphics[width=\textwidth]{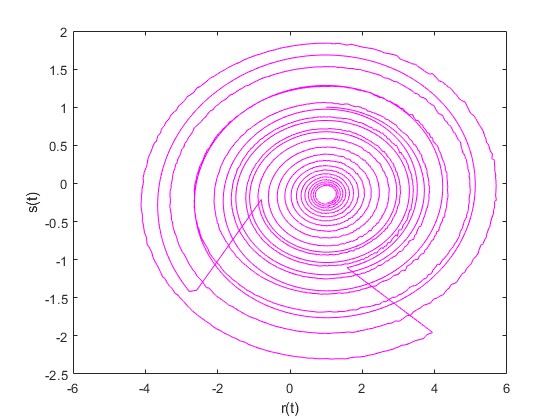}
		\caption{$\alpha=0.99$, $\epsilon=0$.}
	\end{subfigure}
	\hfill
	\begin{subfigure}{0.32\textwidth}
		\includegraphics[width=\textwidth]{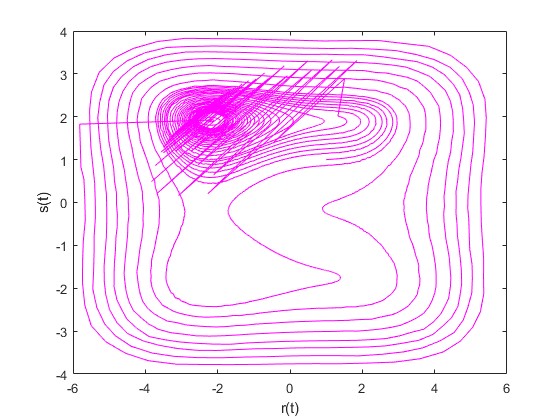}
		\caption{$\alpha=0.94$, $\epsilon=0.3$.}
	\end{subfigure}
	\hfill
	\begin{subfigure}{0.32\textwidth}
		\includegraphics[width=\textwidth]{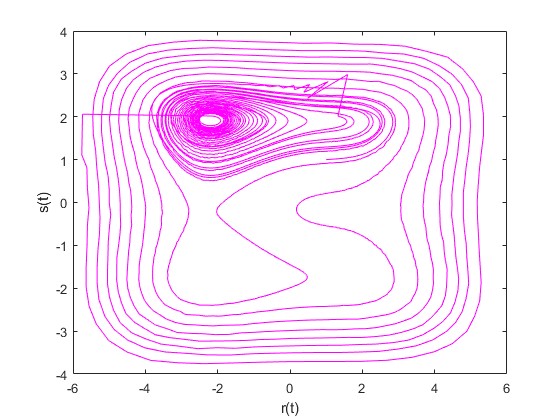}
		\caption{$\alpha=0.98$, $\epsilon=0.3$.}
	\end{subfigure}
	\caption{Chaotic dynamics for the second model \eqref{CFnonlinear} with  $t(0)=0$, $r(0)=1$, $s(0)=1$,  $\rho_{1}=0.12$, $\rho_{2}=0.01$,
		$\phi_{1}=1$, $\phi_{2}=1$,
		$\sigma_{1}=0.01$, $\sigma_{2}=0.02$,
		$\omega_{1}=6.1$, $\omega_{2}=-1$.}
	\label{fig:fig12l}
\end{figure}

It presents the numerical simulations and chaotic phase portrait via first \eqref{CFlinear} and second \eqref{CFnonlinear} models for love dynamical system.
The figures \ref{fig:fig9l}(a)-\ref{fig:fig9l}(l), and \ref{fig:fig10l}(a)-\ref{fig:fig10l}(l) demonstrate the numerical simulations, while the figures \ref{fig:fig11l}(a)-\ref{fig:fig11l}(l), and \ref{fig:fig12l}(a)-\ref{fig:fig12l}(l) demonstrate the chaotic behaviors of phase portrait for the system under consideration.

\section{Conclusion} 
By combining the concepts of piecewise modeling, fractional derivatives, and stochastic analysis, we revised and adjusted the current love dynamical models in order to examine the crossover emotional patterns in romantic love relationships. They have demonstrated an ability to accurately capture both the intricacy of nature and crossover behaviors. In particular, the Atangana-Baleanu derivative, the Caputo derivative, and the Caputo-Fabrizio derivative are the three distinct fractional derivatives that we suggested injecting into the piecewise loving dynamical models. To a specific piecewise loving dynamical model with the Caputo derivative, we proved existence and uniqueness results that may be generalized in a similar way to the other examples. Lagrange quadratic interpolation forms the basis of the numerical approximations of the presented models. Various values of $\alpha$ were simulated numerically for modified models to illustrate.

\end{document}